\documentclass[a4paper,12pt,reqno]{amsart}
\usepackage{amssymb,amsmath,amsfonts,latexsym} 
\usepackage{hyperref}
\usepackage{ifthen}
\usepackage{graphicx}
\usepackage{float}
\usepackage{mathrsfs}
\usepackage{comment}
\usepackage{enumerate}
\usepackage[foot]{amsaddr}
 \usepackage[left=25mm,right=25mm,top=25mm,bottom=25mm,paper=a4paper]{geometry}

\theoremstyle{plain}

\newtheorem*{thm*}{Theorem}
\numberwithin{equation}{section}

\newtheorem{prop}{Proposition}[section]

\theoremstyle{definition}
\newcounter {own}
\def\theown {\thesection  .\arabic{own}}

{\qed\bigskip}

\newcounter{alphabet}

\newcommand{\R}{\mathbb{R}}

\newcommand{\Z}{\mathbb{Z}}
\newcommand{\TR}{\tilde{\rho}}

\newtheorem{theorem}{Theorem}[section]
\newtheorem{definition}{Definition}[section]
\newtheorem{remark}{Remark}[section]
\newtheorem{lemma}{Lemma}[section]

\newcounter{minutes}
\divide\time by 60
\newcounter{hours}
\multiply\time by 60 \addtocounter{minutes}{-\time}

\begin{document}
 
 \title {Hardy-Littlewood maximal, generalized Bessel-Riesz and generalized fractional integral operators in generalized Morrey and $BMO_\phi$ spaces associated with Dunkl operator on the real line}

\thanks{
File:~\jobname .tex,
          printed: \number\year-\number\month-\number\day,
          \thehours.\ifnum\theminutes<10{0}\fi\theminutes}

\author{Sumit Parashar}
 \address{Sumit Parashar, School of Mathematical and Statistical Sciences, Indian Institute of Technology Mandi, Kamand,
 Mandi, HP 175005, India \newline
Email: d22035@students.iitmandi.ac.in }

\author{Saswata Adhikari $^\dagger$}
\address{Saswata Adhikari, School of Mathematical and Statistical Sciences, Indian Institute of Technology Mandi, Kamand,
 Mandi, HP 175005, India. \newline
Email: saswata@iitmandi.ac.in \newline
$^\dagger$ {\tt Corresponding author}}





\begin{abstract}
The analysis of Morrey spaces, generalized Morrey spaces and $BMO_\phi$ spaces related to the Dunkl operators on $\mathbb{R}$ are covered in this paper. We prove the boundedness of the Hardy-Littlewood maximal operators, Bessel-Riesz operators, generalized Bessel-Riesz operators, and generalized fractional integral operators associated with Dunkl operators on $\mathbb{R}$ in the generalized Dunkl-type Morrey spaces. Further, we derive the boundedness of the modified version of the generalized fractional integral operators associated with the Dunkl operators on $\mathbb{R}$ in Dunkl-type $BMO_\phi$ spaces. 
\end{abstract}

\maketitle
\pagestyle{myheadings}
\markboth{Sumit Parashar and Saswata Adhikari}{Hardy-Littlewood Maximal, Generalized Bessel-Riesz, ...}

\subjclass \textbf{Mathematics Subject Classification[2010]} {Primary  42B10, 47G10; Secondary 42B20, 42B25}

\section{Introduction} \label{sec:intro}
The Bessel-Riesz operator $I_{\beta, \gamma}( 0 < \beta < n ~~~ and~~ ~\gamma \geq 0)$ is defined by 
\begin{equation}
    I_{\beta, \gamma} f(x) := \int_{\mathbb{R}^n} f(y) K_{\beta, \gamma}(x -y) dy, \label{eq: 1.1}
\end{equation}
for every $f \in L^p_{loc}(\mathbb{R}^n),~ p \geq 1$, where $K_{\beta, \gamma}(x) := \frac{|x|^{\beta - n}}{(1 + |x|)^\gamma}$. Here $K_{\beta, \gamma}$ is known as Bessel-Riesz kernel. Note that $I_{\beta, 0}$ agrees with the fractional integral operator $I_\beta$ with kernel $K_\beta(x) = |x|^{\beta - n}$.

For $\TR : \mathbb{R}^+ \to \mathbb{R}^+ $ and $\gamma \geq 0$, the generalized Bessel-Riesz operator $I_{\TR, \gamma}$ is defined by
\begin{align*}
    I_{\TR, \gamma} f(x) := \int_{\mathbb{R}^n} \frac{\TR(|x -y|)}{(1 + |x -y|)^\gamma} f(y) dy.
\end{align*}
If $\TR(t) = t^{\beta - n} ~ (0 < \beta < n)$, then $I_{\TR, \gamma}$ reduces to  $I_{\beta, \gamma}$, which is the generalized Bessel-Riesz operator  defined in \eqref{eq: 1.1}.

For a function $\rho : \mathbb{R}^+ \to \mathbb{R}^+$, the generalized fractional integral operator is defined by 
\begin{align*}
    T_\rho f(x) := \int_{\mathbb{R}^n} f(y) \frac{\rho(|x -y|)}{|x -y|^n} dy,
\end{align*}
and its modified version $\tilde{T}_\rho$ is defined by 
\begin{align*}
\tilde{T}_\rho f(x) := \int_{\mathbb{R}^n} f(y) \left( \frac{\rho(|x -y|)}{|x-y|^n} - \frac{\rho(|y|)}{|y|^n}(1 - \chi_{B_0}(y) ) \right) dy,
\end{align*}
where $B_0$ and $\chi_{B_0}$ denote the unit ball around the origin and the characteristic function of $B_0$, respectively. If $\rho(t) = t^\beta , 0 < \beta < n$, then $T_\rho = I_\beta$(fractional integral operator).

For $1 \leq p < \infty$ and a function $\phi : \mathbb{R}^+ \to \mathbb{R}^+$, the generalized Morrey space on $\mathbb{R}^n$ is defined by 
\begin{align*}
    L^{p, \phi}(\mathbb{R}^n) := \{ f \in L^p_{loc}(\mathbb{R}^n) : ||f||_{L^{p, \phi}} < \infty \},
\end{align*}
where $||f||_{L^{p, \phi}} := \displaystyle \sup_{\substack{{r>0}\\{a \in \mathbb{R}^n}}} \frac{1}{\phi(r)} \left ( \frac{1}{r^n} \int_{B(a,r)} |f(x)|^p dx \right)^{\frac{1}{p}} $. If $\phi(r) = r^{- \frac{n}{q}}$, where $ 1 \leq p \leq q < \infty$, then $L^{p, \phi}(\mathbb{R}^n)$ is the Morrey space $L^{p, q}(\mathbb{R}^n)$ and moreover for $p=q$ it is the Lebesgue space $L^p(\mathbb{R}^n)$. 
For a function $\phi : \mathbb{R}^+ \to \mathbb{R}^+$, the $BMO_\phi(\mathbb{R}^n)$ space is defined by 
\begin{align*}
    BMO_\phi(\mathbb{R}^n) := \{ f \in L^1_{loc}(\mathbb{R}^n) : ||f||_{BMO_\phi} < \infty \},
\end{align*}
where $$ ||f||_{BMO_\phi} := \displaystyle \sup_{\substack{{r>0}\\{a \in \mathbb{R}^n}}} \frac{1}{\phi(r)} \frac{1}{r^n} \int_{B(a, r)} | f(x) - f_{B(a,r)}| dx $$  with $f_{B(a,r)} = \frac{1}{r^n} \int_{B(a, r)} f(x) dx$, where $B(a,r)$ denotes the ball center with $a$ and radius $r$. If $\phi(r) \equiv 1$, then $BMO_\phi(\mathbb{R}^n) = BMO(\mathbb{R}^n)$ is the space of bounded mean oscillation functions.

The boundedness of the fractional integral operator $I_\beta$  on Lebesgue spaces was proved by Hardy and Littlewood in \cite{hardy1928some}, \cite{hardy1932some} and by Sobolev in \cite{sobolev1938theorem}. Further, the boundedness of $I_\beta$ on Morrey spaces was studied in \cite{adams1975note}, \cite{chiarenza1987morrey}, \cite{peetre1969theory}. The boundedness of the Bessel-Riesz operator $I_{\beta, \gamma}$ on generalized Morrey spaces was studied by Idris et al. in \cite{idris2016boundedness} and the boundedness of the generalized Bessel-Riesz operator $I_{\TR, \gamma}$ on generalized Morrey spaces was proved by Idris and Gunawam \cite{2017boundedness}. Nakai 
\cite{eridani2002boundedness}, \cite{eridani2004generalized} proved the boundedness of the generalized fractional integral operator $T_{\rho}$ on generalized Morrey spaces. The boundedness of the modified version of the generalized fractional integral operator $\tilde{T}_\rho$ from $BMO_\phi(\mathbb{R}^n)$ to $BMO_\psi(\mathbb{R}^n)$ was proved by Nakai \cite{nakai2001generalized}. Recently, the fractional integral operator on generalized Morrey spaces was studied associated with ball Banach function spaces in \cite{wei}. We also refer \cite{Gul}, \cite{Ruz} for analysis local 
(central) Morrey spaces and generalized local (central) Morrey spaces on homogeneous groups. 

In $1989$, the Dunkl operators were introduced by Dunkl \cite{dunkl1989differential}, which are differential-difference operators on the real line. For a real parameter $\alpha \geq -\frac{1}{2}$, the Dunkl operators are denoted by $\Lambda_\alpha$ and these are associated with the reflection group $Z_2$ on $\mathbb{R}$.  We also refer \cite{de1993dunkl}, \cite{rsl2}, \cite{sy} for more details on Dunkl theory. Using the Dunkl kernel, Dunkl defined the Dunkl transform $\mathcal{F}_\alpha$ in \cite{dunkl1992hankel}. R\"{o}sler in \cite{rosler1994bessel} showed that the Dunkl kernel verifies a product formula. This allows one to define the Dunkl translation $\tau_x^\alpha$, $x \in \mathbb{R}$ and as a result one has the Dunkl convolution.

The maximal operator and the fractional maximal operator were studied associated with Dunkl operator on $\mathbb{R}$ in \cite{abdelkefi2007dunkl}, \cite{guliyev2010p} and \cite{thangavelu2007riesz}, respectively. Guliyev et al. \cite{guliyev2010necessary}  studied the boundedness of the Dunkl-type fractional maximal operator in the Dunkl-type Morrey space. Adhikari and Parui \cite{sps} proved the boundedness of Dunkl-type Bessel-Riesz operators in  Dunkl-type Morrey spaces. Guliyev and Mammadov \cite{guliyev2009fractional} proved that the Dunkl-type modified fractional integral operator is bounded from the space $L^p(\mathbb{R}, d\mu_\alpha)$ to the Dunkl-type $BMO$ space i.e., $BMO(\mathbb{R}, d\mu_\alpha)$.
 
In this paper we are interested in the boundedness of the Bessel-Riesz operator, the generalized Bessel-Riesz operator and the generalized fractional integral operator associated with the Dunkl operator on $\mathbb{R}$ in generalized Dunkl-type Morrey spaces. We also prove that the modified version of the generalized fractional integral operator associated with the Dunkl operator is bounded from $BMO_\phi(\mathbb{R}, d\mu_\alpha)$ to $BMO_\psi(\mathbb{R}, d\mu_\alpha)$.

For the convenience of the reader, let us now shortly summarize the main results of this paper. The following formulations provide the definitions for the spaces that we shall include in this paper, see \eqref{eq: 3.1a} for Dunkl-type Morrey spaces $L^{p, q}(\mathbb{R}, d\mu_\alpha)$, \eqref{eq: 3.2a} for generalized Dunkl-type Morrey spaces $L^{p, \phi}(\mathbb{R}, d\mu_\alpha)$ and \eqref{eq: 7.1} for Dunkl-type $BMO_\phi$ spaces, as well as \eqref{eq: 3.4a} for the Dunkl-type Hardy-Littlewood maximal operator $M^\alpha$, \eqref{eq: 2.1} for Dunkl-type Bessel-Riesz operators $I_{\beta, \gamma}^\alpha$. \eqref{eq:34} for generalized Dunkl-type Bessel-Riesz operators $I_{\TR, \gamma}^\alpha$, \eqref{eq:38} for generalized Dunkl-type fractional integral operators $T_\rho^\alpha$ and \eqref{eq: 7.2} for the modified version of the generalized Dunkl-type fractional integral operators $\tilde{T}_\rho^\alpha$. 

Thus, in this paper we show the following properties:
\begin{itemize}
    \item  Let $1 < p < \infty$ and $\phi(r)$ be a positive measurable function on $\mathbb{R}^+$. Then for $ f \in L^{p, \phi}(\R, d\mu_\alpha)$, we have 
\begin{align}
    ||M^\alpha f||_{L^{p, \phi}(\R, d\mu_\alpha)} \leq C ||f||_{L^{p, \phi}(\R, d\mu_\alpha)}. \nonumber
\end{align} 
\item Let $\gamma > 0$ and $0 < \beta < d_\alpha$. Then $K_{\beta, \gamma}^\alpha \in L^t(\R, d\mu_\alpha)$ and 
        $$ ||K_{\beta, \gamma}^\alpha||_{L^t(\R, d\mu_\alpha)} \sim \left(\sum_{k \in \Z} \frac{(2^kR)^{(\beta - d_\alpha)t + d_\alpha}}{(1 + 2^kR)^{\gamma t}}\right)^{\frac{1}{t}},
        $$
         for  $\frac{d_\alpha}{d_\alpha + \gamma - \beta} < t < \frac{d_\alpha}{d_\alpha - \beta}$ and any $R > 0$.
\item   Let $ 0 < \beta < d_\alpha$ and $\gamma > 0$. If $\phi(r) \leq C r^\nu$ for every $r > 0$, with $\nu < -\beta$, then for $1 < p < \infty$ and $t \in \left( \frac{d_\alpha}{d_\alpha + \gamma - \beta}, \frac{d_\alpha}{d_\alpha - \beta}\right)$, the Dunkl-type Bessel-Riesz operator $I_{\beta, \gamma}^\alpha $ is bounded from  $L^{p, \phi}(\R, d\mu_\alpha)$ to  $L^{q, \psi}(\R, d\mu_\alpha)$ i.e.
    \begin{align*}
        ||I_{\beta, \gamma}^\alpha f||_{L^{q, \psi}(\R, d\mu_\alpha)} \leq C ||K_{\beta, \gamma}^\alpha||_{L^t(\R, d\mu_\alpha)} ||f||_{L^{p, \phi}(\R, d\mu_\alpha)},
    \end{align*}
    where $q = \frac{\nu t' p}{\nu t' + d_\alpha}$ and $\psi(r) = (\phi(r))^{\frac{p}{q}}$.
\item    Let $ 0 < \beta < d_\alpha$ and $\gamma > 0$. If $\phi(r) \leq C r^\nu$ for every $r > 0$, with $\nu < -\beta$, $1 < p < \infty$ and $ \frac{d_\alpha}{d_\alpha + \gamma - \beta} < s \leq t < \frac{d_\alpha}{d_\alpha - \beta}$ and $1 \leq s \leq t$,  then for all $f \in L^{p, \phi}(\R, d\mu_\alpha)$ we get
    \begin{align}
        ||I_{\beta, \gamma}^\alpha f||_{L^{q, \psi}(\R, d\mu_\alpha)} \leq C ||K_{\beta, \gamma}^\alpha||_{L^{s, t}(\R, d\mu_\alpha)} ||f||_{L^{p, \phi}(\R, d\mu_\alpha)}, \nonumber
    \end{align}
    where $q = \frac{\nu t' p}{\nu t' + d_\alpha}$ and $\psi(r) = (\phi(r))^{\frac{p}{q}}$.
\item  Assume that $\omega : \R^{+} \to \R^{+}$ satisfies the doubling condition and $ C' r^{\beta - d_\alpha} \leq \omega(r) \leq C r^{-\beta}$ for every $r > 0$. Let $\gamma > 0$ and $ 0 < \beta < d_\alpha$. Further assume that $\frac{d_\alpha}{d_\alpha + \gamma -\beta} < s < \frac{d_\alpha}{d_\alpha - \beta}$, $s \geq 1$. If $\phi(r) \leq C r^\nu $ for every $r > 0$, where $\nu < - \beta < -d_\alpha - \nu$, then for all $f \in L^{p, \phi}(\R, d\mu_\alpha)$, we have
      \begin{align}
          ||I_{\beta, \gamma}^\alpha f||_{L^{q, \psi}(\R, d\mu_\alpha)} \leq C ||K_{\beta, \gamma}^\alpha||_{L^{s, \omega}(\R, d\mu_\alpha)} ||f||_{L^{p, \phi}(\R, d\mu_\alpha)}, \nonumber
    \end{align}
      where $1 < p < \infty$, $q = \frac{\nu p}{\nu + d_\alpha -\beta }$ and $\psi(r) = (\phi(r))^\frac{p}{q}$. 
\item    Let $\gamma > 0$. Further assume that $\TR$ and $\phi$ satisfy the doubling condition \eqref{eq: 3.3a}. Let $\phi$ be surjective and for $1 < p < q < \infty$, it satisfies
    \begin{align}
        \phi(r) \int_0^r \frac{\TR(t)}{t^{\gamma - d_\alpha + 1}} dt + \int_r^\infty \frac{\TR(t) \phi(t)}{t^{\gamma - d_\alpha +1}} dt \leq C (\phi(r))^{\frac{p}{q}},  \nonumber
    \end{align}
    for all $r > 0$. Then the generalized Dunkl-type Bessel-Riesz operator $I_{\TR, \gamma}^\alpha$ is bounded from $L^{p, \phi}(\mathbb{R}, d\mu_\alpha)$ to $L^{q, \psi}(\mathbb{R}, d\mu_\alpha)$ i.e.
    \begin{align}
        ||I^\alpha_{\TR, \gamma}f ||_{L^{q, \psi}(\R, d\mu_\alpha)} \leq C ||f||_{L^{p, \phi}(\R, d\mu_\alpha)},  \nonumber
    \end{align}
    where $\psi(r) = (\phi(r))^{\frac{p}{q}}$.
\item  Let $\rho$ and $\phi$ satisfy the doubling condition \eqref{eq: 3.3a}. Let $\phi$ be surjective and for $1 < p < q < \infty$, it satisfies
    \begin{align}
        \phi(r) \int_0^r \frac{\rho(t)}{t} dt + \int_r^\infty \frac{\rho(t) \phi(t)}{t} dt \leq C (\phi(r))^{\frac{p}{q}},  \nonumber
    \end{align}
    for all $r > 0$. Then the generalized Dunkl-type fractional integral operator is bounded from $L^{p, \phi}(\mathbb{R}, d\mu_\alpha)$ to $L^{q, \psi}(\mathbb{R}, d\mu_\alpha)$ i.e. 
    \begin{align}
        ||T^\alpha_{\rho}f ||_{L^{q, \psi}(\R, d\mu_\alpha)} \leq C ||f||_{L^{p, \phi}(\R, d\mu_\alpha)},  \nonumber
    \end{align}
    where $\psi(r) = (\phi(r))^{\frac{p}{q}}$.
\item    Let $\rho$ satisfy \eqref{eq: 39}, \eqref{eq: 3.3a}, \eqref{eq: 45}, \eqref{eq: 46}. Let $\phi$ and $\psi$ be almost increasing, $\phi(r) \sim \phi(2r) $ and $\psi(r) \sim \psi(2r)$. If for all $r > 0$,
      \begin{align}
          \int_r^\infty \frac{\rho(t) \phi(t)}{t^2} dt \leq A \frac{\rho(r) \phi(r)}{r}, \nonumber\\
          \int_0^r \frac{\rho(t)}{t} dt ~\phi(r) \leq A' \psi(r), \nonumber
      \end{align}
      then $\tilde{T}_\rho^\alpha$ is bounded from $BMO_\phi(\mathbb{R}, d\mu_\alpha)$ to $BMO_\psi(\mathbb{R}, d\mu_\alpha)$.
\end{itemize}

The paper is arranged as follows. In Section \ref{sec:2}, we briefly recall the concepts of Dunkl theory on the real line and some known results. In Section \ref{sec:3}, we define the generalized Dunkl-type Morrey space and show the boundedness of the Dunkl-type Hardy-Littlewood maximal function in the generalized Dunkl-type Morrey space. In Section \ref{sec:4}, we prove the boundedness of Dunkl-type Bessel-Riesz operators in generalized Dunkl-type Morrey spaces. The boundedness of the generalized Dunkl-type Bessel-Riesz operator in generalized Dunkl-type Morrey spaces is proved in Section \ref{Sec:5}. In Section \ref{Sec:6}, we define the generalized Dunkl-type fractional integral operator and study its boundedness in generalized Dunkl-type Morrey spaces. Finally, in Section \ref{Sec:7}, the boundedness of the modified version of the generalized Dunkl-type fractional integral operator is proved in Dunkl-type $BMO_\phi$ space. Finally, throughout the paper, $C$ denotes a suitable positive constant that need not be same in every place.

 \section{Preliminaries} \label{sec:2}
 For  a fixed real number $\alpha \geq -\frac{1}{2}$, the Dunkl operator associated with $\alpha$ is defined by
 \begin{align*}
     \Lambda_\alpha f(x) = \frac{df}{dx}(x) + \frac{2 \alpha +1}{x} \frac{f(x) - f(-x)}{2},  ~x \in \mathbb{R},
 \end{align*}
 where $f \in C^\infty(\mathbb{R})$. Note that the Dunkl derivative $\Lambda_{-\frac{1}{2}}$ is equal to the classical derivative $\frac{d}{dx}$.
 
 For $\alpha \geq - 1/2$ and $\lambda \in \mathbb{C}$, the initial value problem $\Lambda_\alpha f(x) = \lambda f(x)$; $f(0) =1$ has a unique solution $E_\alpha(\lambda x)$ called Dunkl kernel (\cite{de1993dunkl}, \cite{dunkl1991integral},  \cite{sifi2002generalized}), which is  given by 
 \begin{align*}
     E_\alpha(\lambda x) = \mathcal{J}_\alpha (\lambda x) + \frac{\lambda x}{2(\alpha + 1)} \mathcal{J}_{\alpha+1} (\lambda x), ~ x \in \mathbb{R}, 
 \end{align*}
 where
 \begin{align*}
     \mathcal{J}_\alpha (\lambda x) := \Gamma (\alpha +1) \sum_{n=0}^\infty \frac{(\lambda x / 2)^{2n}}{n! \Gamma(n + \alpha + 1)}
 \end{align*}
 is the modified spherical Bessel function of order $\alpha$. Note that $E_{- 1/ 2}(\lambda x) = e^{\lambda x}$. 
 
 Let $\mu_\alpha$ be the weighted Lebesgue measure on $\mathbb{R}$ given by 
 \begin{align}
     d\mu_\alpha(x) :=  A_\alpha |x|^{2 \alpha + 1} dx,  \nonumber
 \end{align}
 where $A_\alpha = (2 ^ {\alpha + 1 } \Gamma(\alpha + 1))^{-1}$. For $1 \leq p \leq \infty$, the space $L^p(\mathbb{R}, d\mu_\alpha)$ denotes the class of  complex-valued measurable functions $f$ on $\mathbb{R}$ such that 
 \begin{align}
 ||f||_{L^p(\mathbb{R}, d\mu_\alpha)} :=\left(\int_{\mathbb{R}} |f(x)|^p d\mu_\alpha(x)\right)^{\frac{1}{p}} < \infty ,~ ~ if ~ ~ p \in [1, \infty) \nonumber
 \end{align}
 and 
 \begin{align}
 ||f||_{L^{\infty}(\mathbb{R}, d\mu_\alpha)} := ess \displaystyle \sup_{x \in \mathbb{R}} |f(x)| < \infty, ~ ~ if ~ ~ p = \infty. \nonumber
 \end{align}
 
 The Dunkl kernel gives rise to an integral transform, called the Dunkl transform on $\mathbb{R}$, which was introduced in \cite{de1993dunkl}.
 The Dunkl transform of a function $f \in L^1(\mathbb{R}, d\mu_\alpha)$, is given by 
 \begin{align*}
     \mathcal{F}_\alpha (f)(\lambda) := \int_{\mathbb{R}} E_\alpha(- i \lambda x) f(x) d\mu_\alpha(x), ~ \lambda \in \mathbb{R}.
 \end{align*}
 The above integral makes sense as $|E_\lambda (i x)| \leq 1$, $\forall~ x \in \mathbb{R}$ (see \cite{rosler1994bessel}). Note that $\mathcal{F}_{- 1/ 2}$ agrees with the classical Fourier transform $\mathcal{F}$, which is given by
 \begin{align*}
     \mathcal{F}(f)(\lambda) := (2 \pi)^{-1/2} \int_{\mathbb{R}} e^{-i \lambda x} f(x) dx, ~ ~ \lambda \in \mathbb{R}.
 \end{align*}
 
 The Dunkl transform $\mathcal{F}_\alpha$ satisfies the following properties (see \cite{trime2002paley}):
 \begin{itemize}
     \item[(i)] For all $ f \in L^1(\mathbb{R}, d\mu_\alpha)$, one has $|| \mathcal{F}_\alpha (f)||_{L^\infty(\mathbb{R}, d\mu_\alpha)} \leq ||f||_{L^1(\mathbb{R}, d\mu_\alpha)}$.
     \item [(ii)] For all $f \in \mathcal{S}(\R)$,
     \begin{align*}
         \mathcal{F}_\alpha ( \Lambda_\alpha f)(\lambda) = i \lambda \mathcal{F}_\alpha (f)(\lambda), ~~ \lambda \in \mathbb{R}.
     \end{align*}
     \item [(iii)] The Dunkl transform $\mathcal{F}_\alpha$ is an isometric isomorphism from $L^2(\mathbb{R}, d\mu_\alpha)$ onto $L^2(\mathbb{R}, d\mu_\alpha)$. In particular, it satisfies the Plancherel formula i.e.,
     \begin{align*}
         ||\mathcal{F}_\alpha (f)||_{L^2(\mathbb{R}, d\mu_\alpha)} = ||f||_{L^2(\mathbb{R}, d\mu_\alpha)}.
     \end{align*}
     \item[(iv)] If $f \in L^1(\mathbb{R}, d\mu_\alpha)$ with $\mathcal{F}_\alpha(f) \in L^1(\mathbb{R}, d\mu_\alpha)$, then one has the inversion formula 
     \begin{align*}
         f(x) = \int_{\mathbb{R}} \mathcal{F}_\alpha (f) (\lambda) E_\alpha(i \lambda x) d\mu_\alpha(\lambda), ~~ a.e.~~ x \in \mathbb{R}.
     \end{align*}
 \end{itemize}
 
Now we will denote $B(x,R) = \{y \in \mathbb{R} : |y| \in ]max\{0, |x|-R\}, |x| + R[\}$, $R > 0 $, if $x\neq 0$ and $B(0, R) = ]-R, R[ $. Then $\mu_\alpha (B(0, R)) = b_\alpha R^{d_\alpha}$, where $b_\alpha = [2^{\alpha + 1}(\alpha + 1) \Gamma(\alpha + 1)]^{-1}$ and $d_\alpha = 2 \alpha + 2$. If $|x| \leq r$ then $B(x,r) = B(0, |x| + r) \subseteq B(0, 2r)$ and if $|x| > r$ then $B(x,r)= \{ y \in \mathbb{R} : |y| \in ]|x|-r, |x|+r[ \}\subseteq B(0, |x| +r) \subseteq B(0, 2|x|)$.
 
 For $x, y, z \in \mathbb{R}$, we put $$ W_\alpha(x, y, z) = (1 - \sigma_{x, y, z} + \sigma_{z, x, y} + \sigma_{z, y, x} ) \kappa_\alpha(|x|, |y|, |z|)$$ 
 where 
\begin{align}
 \sigma_{x, y, z} &:=
 \begin{cases}
     \frac{x^2 + y^2 - z^2}{2xy}, & if ~ x, y  \in \mathbb{R}\setminus \{0\},\\
     0, & \text{otherwise}, 
 \end{cases}\nonumber
\end{align} 
    and $\kappa_\alpha$ is the Bessel kernel given by 
\begin{align}    
    \kappa_\alpha(|x|, |y|, |z|) :=
    \begin{cases}
        D_\alpha \frac{([(|x| + |y|)^2 - z^2][z^2 - (|x| - |y|)^2])^{\alpha - \frac{1}{2}}}{|xyz|^{2\alpha}}, & if ~ |z| \in S_{x, y},\\
        0, & \text{otherwise},
    \end{cases} \nonumber
\end{align}
where $D_\alpha = (\Gamma(\alpha + 1))^2/ (2^{\alpha - 1} \sqrt{\pi} \Gamma(\alpha + \frac{1}{2}))$ and $ S_{x, y} = [||x|-|y||,|x|+|y|]$.\\ 
Then the signed measure $\nu_{x, y}$ on $\mathbb{R}$ is given by
  \begin{align}
      d\nu_{x, y}(z) &= \begin{cases}
      W_\alpha(x, y, z) d\mu_\alpha(z), & if ~ x, y \in \R \backslash \{0\},\\ d\delta_x(z), & if ~ y = 0,\\ d\delta_y(z), & if ~ x = 0,
  \end{cases} \nonumber
\end{align} 
where the signed kernel  $W_\alpha$ is an even function with respect to all variables satisfying the following properties (see \cite{rosler1994bessel}):
\begin{align}
& W_\alpha(x, y, z) = W_\alpha(y, x, z) = W_\alpha(-x, z, y), \nonumber
\\
& W_\alpha(x, y, z) = W_\alpha(-z, y, -x) = W_\alpha(-x, -y, -z) \nonumber
   \end{align}
   and 
   \begin{align} \int_{\mathbb{R}} |W_\alpha(x, y, z)| d\mu_\alpha(z) \leq 4. \nonumber
   \end{align}
     Let $x, y \in \R$ and $f$ be a continuous function on $\R$. Then the Dunkl translation operator of $f$ is defined by
     \begin{equation}
         \tau_x^\alpha f(y) = \int_{\R} f(z) d\nu_{x,y}(z).  \label{2.1a}
     \end{equation}
Let $z= (x, y)_\theta = \sqrt{x^2 + y^2 - 2 xy cos \theta}$. Then by change of variable we also have (see \cite{abdelkefi2007characterization}) 
    \begin{align}
        \tau_x^\alpha f(y) := \frac{c_\alpha}{2} \int_0^\pi &\bigg\{[f((x,y)_\theta) + f(-(x,y)_ \theta)] \nonumber
        \\
        &\quad +\frac{x+y}{(x,y)_\theta} [f((x,y)_\theta) - f(-(x,y)_\theta)]\bigg\}(1 - cos \theta) sin^{2 \alpha}\theta d\theta, \label{2.1b}
    \end{align}
    where $ c_\alpha =  \left(\int_0^\pi sin^{2 \alpha}\theta d\theta \right)^{-1} =\frac{\Gamma(\alpha + 1)}{\sqrt{\pi}\Gamma(\alpha + \frac{1}{2})}. $
    
     The Dunkl translation operator satisfies the following properties:
      \begin{prop} \cite{mourou2001transmutation} \label{pro: 2.1}
        \begin{itemize}
            \item [(i)] $\tau_x^\alpha, ~x \in \mathbb{R}$ is a bounded linear operator on $C^\infty(\mathbb{R})$
            \item [(ii)] For all $f \in C^\infty(\mathbb{R})$ and $x,y \in \mathbb{R}$, one has 
            \begin{align*}
                \tau_x^\alpha f(y) = \tau_y^\alpha f(x), \tau_0^\alpha f(x) = f(x), ~ \tau_x^\alpha \circ \tau_y^\alpha = \tau_y^\alpha \circ \tau_x^\alpha.
            \end{align*} 
        \end{itemize}
      \end{prop}

 If $f, g \in L^2(\mathbb{R}, d\mu_\alpha)$, then 
            \begin{align}
                \int_{\mathbb{R}} \tau_x^\alpha f(y) g(y) d\mu_\alpha(y) = \int_{\mathbb{R}} f(y) \tau_{-x}^\alpha g(y) d\mu_\alpha(y), ~ \forall~ x \in \mathbb{R}. \nonumber
            \end{align}
            
Let $f$, $g$ are continuous functions on $\mathbb{R}$ with compact support. The generalized convolution $f*_\alpha g$ is defined by 
\begin{align*}
    f*_\alpha g(x) := \int_{\mathbb{R}} f(y) \tau_y^\alpha g(x) d\mu_\alpha(y). 
\end{align*}
The generalized convolution $*_\alpha$ is associative and commutative \cite{rosler1994bessel}. 

We also mention in the following some useful properties of the generalized convolution $*_\alpha$.
\begin{theorem}\cite{soltani2004lp} \label{tm: 2.1}
\begin{enumerate}
    \item[(i)]   For all $x \in \mathbb{R}$ and $p \geq 1$, the generalized translation operator $\tau_x^\alpha$ is bounded from $L^p(\mathbb{R}, d\mu_\alpha)$ to $L^p(\mathbb{R}, d\mu_\alpha)$ i.e.
      \begin{align*}
          ||\tau_x^\alpha f||_{L^p(\mathbb{R}, d\mu_\alpha)} \leq 4 ||f||_{L^p(\mathbb{R}, d\mu_\alpha)}.
       \end{align*}
    \item [(ii)] If $f\in L^{1}(\mathbb{R}, d\mu_\alpha)$, then 
\begin{eqnarray*}
\mathcal{F}_{\alpha} (\tau_x^\alpha f)(\lambda)=E_{\alpha}(i\lambda x)\mathcal{F}_{\alpha}(f)(\lambda),~x,\lambda\in\mathbb{R}.
\end{eqnarray*}
\item [(iii)] If $f\in L^{1}(\mathbb{R}, d\mu_\alpha)$ and $g\in L^{2}(\mathbb{R}, d\mu_\alpha)$, then 
\begin{eqnarray*}
\mathcal{F}_{\alpha}(f\ast_\alpha g)=\mathcal {F}_{\alpha} (f)\mathcal{F}_{\alpha} (g).
\end{eqnarray*}
 \item[(iv)] Assume that $p,q,r \in [1, \infty]$ satisfying $\frac{1}{p} + \frac{1}{q} = 1 + \frac{1}{r}$  (the Young condition ). Then the map $(f,g) \to f *_\alpha g$ defined on $C_c(\mathbb{R}) \times C_c(\mathbb{R})$, extends to a continuous map from
 $L^p(\mathbb{R}, d\mu_\alpha) \times L^q(\mathbb{R}, d\mu_\alpha)$ to $L^r(\mathbb{R}, d\mu_\alpha)$ and one has
\begin{align*}
    ||f *_\alpha g ||_{L^r(\mathbb{R}, d\mu_\alpha)} \leq 4 ||f||_{L^p(\mathbb{R}, d\mu_\alpha)} ||g||_{L^q(\mathbb{R}, d\mu_\alpha)}.
\end{align*} 
\end{enumerate}
\end{theorem}

Now we will recall a few definitions and the known results in the following, which will be useful in our paper.
\begin{definition} \cite{sps}
  For $1\leq p \leq q < \infty$, the Dunkl-type Morrey space $L^{p,q}(\R, d\mu_\alpha)$ is defined to be the class of all locally $p$-integrable functions $f$ on $\mathbb{R}$ such that \begin{align}
    ||f||_{L^{p,q}(\R, d\mu_\alpha)} :=  \sup_{\substack{{r>0}\\{x \in \R}}} r^{d_\alpha \left(\frac{1}{q}- \frac{1}{p} \right)} \left(\int_{B(0,r)} \tau_x^\alpha |f|^p(y) d\mu_\alpha(y)\right)^{\frac{1}{p}} < \infty. \label{eq: 3.1a}
\end{align}  
\end{definition}

Let $0 < \beta < d_\alpha$, where $d_\alpha = 2 \alpha + 2$ and $\gamma \geq 0$.  For every locally p-integrable functions $f$ on $\mathbb{R}$, the Dunkl-type Bessel-Riesz operator $I_{\beta, \gamma}^\alpha$ is defined by
\begin{align}
    I_{\beta, \gamma}^\alpha f(x) := \int_\R f(y) \tau_x^\alpha K_{\beta, \gamma}^\alpha(y) d\mu_\alpha(y),\label{eq: 2.1}
\end{align}
 where $1 \leq p < \infty$ and $K_{\beta, \gamma}^{\alpha}(x) = \frac{|x|^{\beta - d_\alpha}}{(1 + |x|)^\gamma}$. 
    The kernel $K_{\beta, \gamma}^\alpha $ is called the Dunkl-type Bessel-Riesz kernel. For $\gamma = 0$, the $K_{\beta, 0}^\alpha = K_{\beta}^\alpha$ is the Dunkl-type Riesz kernel and in this case the Dunkl-type Bessel-Riesz operator reduces to Dunkl-type fractional integral operator $I_\beta^\alpha$. 
     \begin{lemma} \cite{lk}\label{7.1}
      If $f \in L^1(\mathbb{R}, d\mu_\alpha)$ and $g \in L^p(\mathbb{R}, d\mu_\alpha)$, $1 \leq p < \infty$,
then
\begin{align}
\tau_{x_0}^\alpha(f *_\alpha g) = \tau_{x_0}^\alpha f *_\alpha g = f *_\alpha \tau_{x_0}^\alpha g,  ~~x_0 \in \mathbb{R}. \nonumber
\end{align}
    \end{lemma}
  \begin{lemma} \cite{idris2016boundedness} \label{lem: 2.1}
    If $b > a > 0$, then $\displaystyle\sum_{k \in \Z} \frac{(u^k)^a}{(1+ u^kR)^b} < \infty$, for every $ u > 1$ and $ R > 0$.    \end{lemma}
    \begin{lemma} (\cite{abdelkefi2007dunkl}, \cite{guliyev2010p}) \label{lem: 2.2}
        Support of $\tau_x^\alpha \chi_{B(0,r)}$ is contained in the ball $B(x,r)$  for all $x \in \R$. Moreover, 
        \begin{align}
            0 \leq \tau_x^\alpha \chi_{B(0, r)}(y) \leq min \left\{ 1, \frac{2 c_\alpha}{2\alpha + 1} \left( \frac{r}{|x|} \right)^{2 \alpha + 1} \right \} , ~\forall~ y \in B(x, r), \nonumber
        \end{align}
        where $c_\alpha = \frac{ \Gamma(\alpha + 1)}{\sqrt{\pi} \Gamma \left(\alpha + \frac{1}{2} \right)}$. In fact, when $|x| \leq r$, one has $0 \leq \tau_x^\alpha \chi_{B(0,r)}(y) \leq 1$ and when $|x| > r$, one has $0 \leq \tau_x^\alpha \chi_{B(0,r)}(y) \leq \frac{2 c_\alpha}{2 \alpha +1} \left(\frac{r}{|x|} \right)^{2 \alpha + 1}$.
    \end{lemma}

    The following lemma was proved in \cite{dg} on $\mathbb{R}^n$, we will state the same result on $\mathbb{R}$.
\begin{lemma} \label{lem:2.222}
  The kernel $\Phi(x,y) = \tau_x^\alpha K_\beta^\alpha(y)$ satisfies the following properties:
  \begin{itemize}
      \item [(i)] $\Phi(x,y) = \Phi(y,x)$.
      \item [(ii)] $\Phi(rx, ty) = r^{\beta - d_\alpha} \Phi(x, \frac{ty}{r})$.
      \item [(ii)] $\Phi(x,y) = \int_{\mathbb{R}} (x^2 + y^2 - 2x\eta)^{\frac{\beta-d_\alpha}{2}} d\mu_y^\alpha(\eta)$,
  \end{itemize}
  where for each $y$ in $\mathbb{R}$, $d\mu_y^\alpha$ is a probability measure on $\mathbb{R}$, whose support is contained in $[-y, y]$.
  
  Note that $\frac{1}{(1+|x|)^\gamma}$ is a decreasing function on $\mathbb{R}$ with no singular point.
  Following the same procedure as in the proof for the integral representation of the kernel $\tau_x^\alpha K_\beta^\alpha(y)$ in \cite{dg},  the kernel $\tau_x^\alpha K_{\beta, \gamma}^\alpha(y)$ has the following integral representation:
  \begin{align}
      \tau_x^\alpha K_{\beta, \gamma}^\alpha(y) = \int_{\mathbb{R}} \frac{(x^2 + y^2 - 2x \eta)^{\frac{\beta - d_\alpha}{2}}}{(1 + \sqrt{x^2 + y^2 - 2x\eta})^\gamma} d\mu_y^\alpha(\eta). \label{eq: 2.2}
  \end{align}
\end{lemma}
    \begin{lemma} \cite{sps} \label{lem: 2.4}
    For $1 \leq s \leq t$ and $\frac{d_\alpha}{d_\alpha + \gamma - \beta} < t < \frac{d_\alpha}{d_\alpha - \beta}$, if $\gamma > 0$ and $0 < \beta < d_\alpha$, the following holds
    \begin{align*}
        ||K_{\beta, \gamma}^\alpha||_{L^{s,t}(\R, d\mu_\alpha)} \leq ||K_{\beta, \gamma}^\alpha||_{L^t(\R, d\mu_\alpha)}. 
    \end{align*}
\end{lemma}

Now for locally integrable functions $f$ on $\mathbb{R}$, the Dunkl-type Hardy-\\
Littlewood maximal operator $M^\alpha$ (see \cite{guliyev2009fractional}) is defined by
\begin{align}
    M^\alpha f(x) := \sup_{r > 0} \frac{1}{\mu_\alpha(B(0,r))} \int_{B(0,r)} \tau_x^\alpha|f|(y) d\mu_\alpha(y), ~ 
   x \in \mathbb{R}.  \label{eq: 3.4a}
\end{align}
\begin{theorem} (\cite{abdelkefi2007dunkl}, \cite{guliyev2009fractional}, \cite{solatani2005})\label{tm: 2.2}
    \begin{itemize}
        \item[(i)] Let $s > 0$ be a real number. Then for all $ f \in L^1(\mathbb{R}, d\mu_\alpha) $, one has 
        \begin{align*}
            \mu_\alpha\{ x \in \mathbb{R} : M^\alpha f(x) > s \} \leq \frac{A_1}{s} \int_{\mathbb{R}} |f(x)| d\mu_\alpha(x),
        \end{align*}
       where $A_1 > 0$ is independent of $f$.
    \item[(ii)]  If $f \in L^p(\mathbb{R}, d\mu_\alpha), 1 < p \leq \infty $, then $M^\alpha f \in L^p(\mathbb{R}, d\mu_\alpha)$ and 
    \begin{align*}
        ||M^\alpha f ||_{L^p(\mathbb{R}, d\mu_\alpha)} \leq A_2 ||f||_{L^p(\mathbb{R}, d\mu_\alpha)},
    \end{align*}
    where $A_2 > 0$ is independent of $f$.
    \end{itemize}
\end{theorem}

 For functions $ \theta_1, \theta_2 :(0,+\infty) \to (0,+\infty)$, we denote $\theta_1(r) \sim \theta_2(r)$ if there exists a constant $C>0$ such that $ \frac{1}{C}\theta_1(r) \leq \theta_2(r) \leq C \theta_1(r)$, $r>0$.\\
 A function $\theta : (0,+\infty) \to (0,+\infty)$ is said to be almost increasing (almost decreasing) if there exists a constant $C>0$ such that $\theta(r) \leq C \theta(s) ~ (\theta(r) \geq C\theta(s))$  for $r \leq s$.
\section{The boundedness of the Dunkl-type Hardy-Littlewood Maximal Operator in generalized Dunkl-type Morrey spaces}\label{sec:3}

First, let us define the generalized Dunkl-type Morrey space on real line. After that we prove that the Dunkl-type Hardy-Littlewood maximal operator is bounded in this space. 

Let $\phi : \R^{+} \to \R^{+}$ be a function which is measurable. For $1\leq p < \infty$, we define the generalized Dunkl-type Morrey space $L^{p, \phi}(\R, d\mu_\alpha)$ for all locally $p$-integrable functions $f$ on $\mathbb{R}$ such that 
\begin{align}
    ||f||_{L^{p, \phi}(\R, d\mu_\alpha)} :=  \sup_{\substack{{r > 0}\\{x \in \R}}} \frac{1}{\phi(r)}\left(\frac{1}{r^{d_\alpha}}\int_{B(0,r)} \tau_x^\alpha|f|^p(y) d\mu_\alpha(y) \right)^{\frac{1}{p}} < \infty. \label{eq: 3.2a}
\end{align} 
We also assume that $\phi$ is almost decreasing and $r^{\frac{d_\alpha}{p}}\phi(r)$ is almost increasing, so that there exists a constant $C_1 > 0$ such that 
\begin{align}
    \frac{1}{2} \leq \frac{r}{s} \leq 2 \implies \frac{1}{C_1} \leq \frac{\phi(r)}{\phi(s)} \leq C_1, \label{eq: 3.3a}
\end{align} 
i.e. $\phi$ satisfies the doubling condition.
In particular, when $\phi(r) = r^{-\frac{d_\alpha}{q}}$, $1 \leq p \leq q < \infty$, then it can be easily verified that the generalized Dunkl-type Morrey space $L^{p, \phi}(\R, d\mu_\alpha)$ reduces to the Dunkl-type Morrey space $L^{p,q}(\R, d\mu_\alpha)$, which is defined in \eqref{eq: 3.1a} and therefore it justifies the name generalized Dunkl-type Morrey space given to $L^{p, \phi}(\R, d\mu_\alpha)$. 
\begin{theorem} \label{tm: 3.1a}
For $1 < p < \infty$, let $\phi$ be a positive measurable function on $\mathbb{R}^+$. 
Then for $ f \in L^{p, \phi}(\R, d\mu_\alpha)$, we have 
\begin{align}
    ||M^\alpha f||_{L^{p, \phi}(\R, d\mu_\alpha)} \leq C ||f||_{L^{p, \phi}(\R, d\mu_\alpha)}.\label{eq: 3.6a}
\end{align}
\end{theorem}
\textbf{Proof.} Since we know that by definition of generalized Dunkl-type Morrey space 
\begin{align}
    ||f||_{L^{p, \phi}(\R, d\mu_\alpha)} &=  \sup_{\substack{{r>0}\\{x \in \R}}} \frac{1}{\phi(r)} \left(\frac{1}{r^{d_\alpha}} \int_{B(0, r)} \tau_x^{\alpha}|f|^p(y) d\mu_\alpha(y) \right)^{\frac{1}{p}}, \nonumber
\end{align}
we have for every $r > 0$ and for every $x \in \mathbb{R}$,
\begin{align}
    \left(\int_{B(0, r)} \tau_x^{\alpha}|f|^p(y) d\mu_\alpha(y) \right)^{\frac{1}{p}} \leq \phi(r) r^{\frac{d_\alpha}{p}} ||f||_{L^{p, \phi}(\R, d\mu_\alpha)}. \label{eq:b}
\end{align}
Now for fixed $r > 0$ and $x_0 \in \R$, we write $f := f\chi_{B(x_0,2r)} + f\chi_{B^c(x_0,2r)} $, where $B^c(x_0,2r)$ is the complement of $B(x_0,2r)$. Then the fact that $M^\alpha$ is sublinear operator leads to 
\begin{align}
   |M^\alpha f|^p(x) &= |M^\alpha (f\chi_{B(x_0,2r)} +f\chi_{B^c(x_0,2r)})|^p(x) \nonumber \\
                     & \leq 2^{p-1} ( |M^\alpha f\chi_{B(x_0,2r)}|^p(x) + |M^\alpha f\chi_{B^c(x_0,2r)}|^p(x)). \nonumber
\end{align}
Therefore, 
\begin{align}
      &\int_{B(0,r)} \tau_{x_0}^\alpha |M^\alpha f|^p(x) d\mu_\alpha(x) \nonumber
      \\
      &= \int_{\mathbb{R}} |M^\alpha f|^p(x) \tau_{-x_0}^\alpha \chi_{B(0,r)}(x) d\mu_\alpha(x) \nonumber
      \\
    &\leq 2^{p-1} \left(\int_{\mathbb{R}} |M^\alpha f\chi_{B(x_0,2r)}|^p(x) \tau_{-x_0}^\alpha \chi_{B(0,r)}(x) d\mu_\alpha(x) \right. \nonumber
   \\
    &\hspace{1cm} \left .+  \int_{\mathbb{R}} |M^\alpha f\chi_{B^c(x_0,2r)}|^p(x) \tau_{-x_0}^\alpha \chi_{B(0,r)}(x) d\mu_\alpha(x) \right). \label{eq: 3.9a}
\end{align}
Making use of Lemma \ref{lem: 2.2} and Theorem \ref{tm: 2.2}, the first integral of the right hand side of \eqref{eq: 3.9a} reduces to
\begin{align}
&\int_\mathbb{R} |M^\alpha f\chi_{B(x_0,2r)}|^p(x) \tau_{-x_0}^\alpha \chi_{B(0,r)}(x) d\mu_\alpha(x) \nonumber
\\
&\leq \int_\mathbb{R} |M^\alpha f\chi_{B(x_0,2r)}|^p(x) d\mu_\alpha(x) \nonumber\\
& \leq C\int_\mathbb{R} |f\chi_{B(x_0,2r)}|^p(x) d\mu_\alpha(x) \nonumber\\
& = C\int_{B(x_0,2r)} |f|^p(x) d\mu_\alpha(x) \nonumber
\\
&= C \int_{I(-x_0, 2r) \cup I(x_0,2r)} |f|^p(x) d\mu_\alpha(x) \nonumber
\\
& \leq C \left(\int_{I(-x_0, 2r)} |f|^p(x) d\mu_\alpha(x) + \int_{I(x_0,2r)} |f|^p(x) d\mu_\alpha(x) \right),\label{eq: 3.10ad}
\end{align}
since $B(x_0, 2r) = I(-x_0, 2r) \cup I(x_0, 2r)$, where $I(x_0, 2r) =(x_0 - 2r, x_0 + 2r)$. Now using the following result (see \cite{yy}, \cite{Nagacy}) 
\begin{align}
    \int_{I(y, \rho)} |f|^p(x) d\mu_\alpha(x) \leq C\int_{B(0,\rho)} \tau_{y}^\alpha |f|^p(x) d\mu_\alpha(x), \nonumber
\end{align}
for all $y \in \R$ and for all $\rho > 0$,
we get from \eqref{eq: 3.10ad}
\begin{align}
 \int_\mathbb{R} |M^\alpha f\chi_{B(x_0,2r)}|^p(x) \tau_{-x_0}^\alpha \chi_{B(0,r)}(x) d\mu_\alpha(x)& \leq C\int_{B(0,2r)} \tau_{x_0}^\alpha|f|^p(x) d\mu_\alpha(x)
\nonumber \\
&\leq C \phi^p(2r) (2r)^{d_\alpha} ||f||^p_{L^{p, \phi}(\mathbb{R}, d\mu_\alpha)} \nonumber 
\\
&\leq C \phi^p(r) r^{d_\alpha} ||f||^p_{L^{p, \phi}(\mathbb{R}, d\mu_\alpha)}.\label{eq: 3.10a}
\end{align}
Now we will find an estimate for the second integral of the right hand side of \eqref{eq: 3.9a}. Towards this we observe that $x \in B(-x_0,r)$. We estimate $M^\alpha$ applied to $f\chi_{B^c(x_0,2r)}$ as
\begin{align}
   &\frac{1}{\phi(r)} M^\alpha f\chi_{B^c(x_0,2r)}(x)\nonumber
   \\
   &= \frac{1}{\phi(r)} \sup_{R > 0} \frac{1}{\mu_\alpha(B(0,R))} \int_{B(0,R)} \tau_x^\alpha |f\chi_{B^c(x_0,2r)}|(y) d\mu_\alpha(y). \nonumber 
\end{align}
 By using support property of $ \tau_{-x}^\alpha\chi_{B(0,R)}$ from the Lemma \ref{lem: 2.2} 
 \begin{align}
      &\frac{1}{\phi(r)} M^\alpha f\chi_{B^c(x_0,2r)}(x) \nonumber
      \\
      &= \frac{1}{\phi(r)} \sup_{R>0} \frac{1}{\mu_\alpha(B(0,R))} \int_\R |f\chi_{B^c(x_0,2r)}|(y) \tau_{-x}^\alpha \chi_{B(0,R)}(y) d\mu_\alpha(y) \nonumber
      \\
      & = \frac{1}{\phi(r)} \sup_{R>0} \frac{1}{\mu_\alpha(B(0,R))} \int_{B(-x,R)\setminus B(x_0,2r)} |f(y)| \tau_{-x}^\alpha \chi_{B(0,R)}(y) d\mu_\alpha(y). \label{eq:max2}
 \end{align}
 Now if $B(-x, R)\setminus B(x_0,2r)$ is an empty set then $M^\alpha f\chi_{B^c(x_0,2r)}(x) = 0$. So the second integral of the right hand side of \eqref{eq: 3.9a} becomes zero. Therefore, with out loss of generality we assume that 
 $B(-x, R)\setminus B(x_0, 2r) \neq \emptyset$. Then we claim that $R > r$. Now we observe that since $B(-x_0,r)= \{ x \in \R : |x| \in ] \max\{ 0, |x_0| -r \}, |x_0| + r[ \}$, we can write
 \begin{align}
     B(-x, R) &= \{ y \in \R : |y| \in ] \max\{ 0, |x| - R\}, |x| + R[\} \nonumber
     \\
     & \subseteq \{ y \in \R : |y| \in ] \max\{ 0, |x_0|-r - R\}, |x_0|+r+ R[\}.
 \end{align}
Therefore, we have  
 \begin{align}
 B(-x, R) \setminus B(x_0, 2r)  \subseteq &\{ y \in \R : |y| \in ] \max\{ 0, |x_0|-r - R\}, |x_0|+r+ R[\} \nonumber
 \\
 &\hspace{0.5cm}\setminus \{ y \in \R : |y| \in ] \max\{ 0, |x_0| - 2r\}, |x_0| + 2r[\}. \label{eq:max3}
 \end{align}
Now we consider the following cases under our assumption $B(-x, R) \setminus B(x_0, 2r) \\
\neq \emptyset$ :
 \begin{itemize}
     \item [(i)] If $\max\{0, |x_0| -2r\}= 0$ and $\max\{0, |x_0| -r -R\}= 0$, then from \eqref{eq:max3}, $B(-x, R) \setminus B(x_0, 2r) \subseteq B(0, |x_0| + r + R) \setminus B(0, |x_0| + 2r) \neq \emptyset$, which implies $|x_0| + r + R > |x_0| + 2r$ i.e. $R > r$.
     \item [(ii)] If $\max\{0, |x_0| -2r\}= 0$ and $\max\{0, |x_0| -r -R\}=  |x_0| - r - R$, then again from \eqref{eq:max3}, $B(-x, R) \setminus B(x_0, 2r) \subseteq \{ y \in \R : |y| \in ] |x_0|-r - R, |x_0|+r+ R[\} \setminus \{ y \in \R : |y| \in ] 0, |x_0| + 2r[\} \subseteq B(0, |x_0| + r + R) \setminus B(0, |x_0| + 2r) \neq \emptyset$, which implies  again $R > r$.
     \item[(iii)]  If $\max\{ 0, |x_0| - 2r \}= |x_0| -2r$ and $\max\{0, |x_0| -r -R\}= 0$, then it follows that $ R +  r > |x_0|$ and $|x_0| > 2r$, which again gives $R > r$.
      \item [(iv)] If $\max\{0, |x_0| - 2r\}= |x_0| - 2r$ and $\max\{0, |x_0| -r -R\}= |x_0| - r-R$, then again from \eqref{eq:max3}, we obtain 
     \begin{align}
         B(-x, R) \setminus B(x_0, 2r) \subseteq & \{ y \in \R : |y| \in ]|x_0| - r - R, |x_0| + r + R[\} \nonumber
         \\
         &\hspace{0.5cm} \setminus \{ y \in \R : |y| \in ]|x_0|-2r, |x_0| + 2r[ \} \neq \emptyset. \nonumber
     \end{align}
    This is possible only if either $|x_0| - r- R < |x_0| -2r$ or $|x_0| + r + R > |x_0| + 2r$, which implies $ R > r$.
 \end{itemize}
 Therefore in all the above cases, we obtain $ R > r$, which proves our claim. Then by virtue of the fact that $\phi$ is almost decreasing, we obtain from \eqref{eq:max2}
\begin{align}
  &\frac{1}{\phi(r)} M^\alpha f\chi_{B^c(x_0,2r)}(x) \nonumber
  \\
  &\leq C\sup_{R > 0} \frac{1}{\phi(R)} \frac{1}{\mu_\alpha(B(0,R))} \int_\R |f(y)| \tau_{-x}^\alpha \chi_{B(0,R)}(y) d\mu_\alpha(y) \nonumber
  \\
  &= C \sup_{R > 0} \frac{1}{\phi(R)} \frac{1}{\mu_\alpha(B(0,R))} \int_\R |f(y)| (\tau_{-x}^\alpha \chi_{B(0,R)}(y))^{\frac{1}{p}+ \frac{1}{q}} d\mu_\alpha(y), \nonumber 
\end{align}
where $q$ is the conjugate exponent of $p \in (1, \infty)$, i.e. $1/p + 1/q = 1$. Applying H\"{o}lder inequality and Theorem \ref{tm: 2.1} , we get
\begin{align}
    &\frac{1}{\phi(r)} M^\alpha f\chi_{B^c(x_0,2r)}(x) \nonumber
    \\
    &\leq C\sup_{R > 0} \frac{1}{\phi(R)} \frac{1}{R^{d_\alpha}} \left(\int_\R |f(y)|^p \tau_{-x}^\alpha \chi_{B(0,R)}(y) d\mu_\alpha(y) \right)^{1/p} \nonumber
    \\
    & \hspace{3.65cm}\times\left(\int_\R \tau_{-x}^\alpha \chi_{B(0,R)}(y) d\mu_\alpha(y) \right)^{1/q}  \nonumber
    \\
    & \leq C\sup_{\substack{{R > 0}\\{x \in \R}}} \frac{1}{\phi(R)} \left( \frac{1}{R^{d_\alpha}}\int_\R |f(y)|^p \tau_{-x}^\alpha \chi_{B(0,R)}(y) d\mu_\alpha(y) \right)^{1/p}  \nonumber
    \\
     & = C\sup_{\substack{{R > 0}\\{x \in \R}}} \frac{1}{\phi(R)} \left( \frac{1}{R^{d_\alpha}}\int_{B(0,R)} \tau_x^\alpha|f(y)|^p d\mu_\alpha(y) \right)^{1/p}  \nonumber
    \\
    & = C ||f||_{L^{p,\phi}(\R, d\mu_\alpha)}. \label{eq: 3.8b}
\end{align}
Now, using \eqref{eq: 3.8b}, we obtain
\begin{align}
     \int_\mathbb{R} |M^\alpha f\chi_{B^c(x_0,2r)}|^p(x) \tau_{-x_0}^\alpha \chi_{B(0,r)}(x) d\mu_\alpha(x) \leq C ||f||^p_{L^{p, \phi}(\mathbb{R}, d\mu_\alpha)} r^{d_\alpha} \phi^p(r). \label{eq: 3.11a}
\end{align}
Substituting the values from \eqref{eq: 3.10a} and \eqref{eq: 3.11a} in \eqref{eq: 3.9a}, we get
\begin{align}
    \int_{B(0,r)} \tau_{x_0}^\alpha |M^\alpha f|^p(x) d\mu_\alpha(x) \leq C ||f||^p_{L^{p, \phi}(\mathbb{R}, d\mu_\alpha)} r^{d_\alpha} \phi^p(r), \nonumber
\end{align}
from which it follows that
 \begin{align}
   \frac{1}{\phi(r)} \left( \frac{1}{r^{d_\alpha}} \int_{B(0,r)} \tau_{x_0}^\alpha |M^\alpha f|^p(x) d\mu_\alpha(x) \right)^{1/p} \leq C ||f||_{L^{p,\phi}(\mathbb{R}, d\mu_\alpha)}. \nonumber
 \end{align}
 Taking the supremum over $r >0$ and $x_0 \in \mathbb{R}$, we have
 \begin{align}
   ||M^\alpha f||_{L^{p, \phi}(\R, d\mu_\alpha)} \leq C ||f||_{L^{p, \phi}(\R, d\mu_\alpha)} ,  \nonumber
 \end{align}
which completes the proof. \qed
\begin{remark}
    The proof of Theorem \Ref {tm: 3.1a} is different in nature because it does not follow in a similar way to that of the Euclidean setting. 
\end{remark}


\section{Inequalities for Dunkl-type Bessel-Riesz operators on generalized Dunkl-type Morrey spaces} \label{sec:4}
In this section, we prove that Dunkl-type Bessel-Riesz operators is bounded on the generalized Dunkl-type Morrey space \eqref{eq: 3.2a}. Towards this, first we prove the following Lemma.
\begin{lemma} \label{lem: 4.1}
        Let $\gamma > 0$ and $0 < \beta < d_\alpha$. Then $K_{\beta, \gamma}^\alpha \in L^t(\R, d\mu_\alpha)$ and 
        $$ ||K_{\beta, \gamma}^\alpha||_{L^t(\R, d\mu_\alpha)} \sim \left(\sum_{k \in \Z} \frac{(2^kR)^{(\beta - d_\alpha)t + d_\alpha}}{(1 + 2^kR)^{\gamma t}}\right)^{\frac{1}{t}},
        $$
         for  $\frac{d_\alpha}{d_\alpha + \gamma - \beta} < t < \frac{d_\alpha}{d_\alpha - \beta}$ and any $R > 0$.
         \end{lemma}
    \textbf{Proof.} For given $R > 0$, we have
\begin{align}
 \int_\R |K_{\beta, \gamma}^\alpha(y)|^t d\mu_\alpha(y) &= A_\alpha \int_\R \frac{|y|^{t (\beta - d_\alpha)}}{(1 + |y|)^{\gamma t}} |y|^{2\alpha + 1} dy, \nonumber
 \\
& =  2 A_\alpha \sum_{k \in \Z} \int_{{2^k R} \leq y < {2^{k+1} R}} \frac{y^{t (\beta - d_\alpha)+ d_\alpha -1}}{(1 + y)^{\gamma t}} dy  \label{eq: 3.1}
 \\
 &\leq   2 A_\alpha \sum_{k \in \Z} \frac{1}{(1 + 2^k R)^{\gamma t}} \int_{{2^k R} \leq y < {2^{k+1} R}} y^{t (\beta - d_\alpha)+d_\alpha -1} dy \nonumber
 \\
 &= B_\alpha \sum_{k \in \Z} \frac{{(2^k R)}^ {t(\beta - d_\alpha)+ d_\alpha}} {(1 + 2^k R)^{\gamma t}} .\nonumber
 \end{align}
 The assumption $\frac{d_\alpha}{d_\alpha + \gamma - \beta} < t < \frac{d_\alpha}{d_\alpha - \beta}$  where $ 0 < \gamma, 0 < \beta < d_\alpha $, implies that $(\beta - d_\alpha) t + d_\alpha > 0$. In Lemma \ref{lem: 2.1}, taking $u = 2$, $a = (\beta - d_\alpha) t + d_\alpha$, $b = \gamma t,$ we get  $\displaystyle\sum_{k \in \Z} \frac{(2^kR)^{(\beta - d_\alpha)t + d_\alpha}}{(1 + 2^kR)^{\gamma t}} < \infty $,  which implies $K_{\beta, \gamma}^\alpha \in L^t(\R, d\mu_\alpha)$. \\
On the other hand, from equation \eqref{eq: 3.1}
\begin{align}
 \int_\R |K_{\beta, \gamma}^\alpha(y)|^t d\mu_\alpha(y) &=  2 A_\alpha \sum_{k \in \Z} \int_{{2^k R} \leq y < {2^{k+1} R}} \frac{y^{t (\beta - d_\alpha) + d_\alpha -1}}{(1 + y)^{\gamma t}} dy \nonumber 
 \\ 
  &\geq \frac{2 A_\alpha}{2 ^ {\gamma t}} \sum_{k \in \Z} \frac{1}{(1 + 2^k R)^{\gamma t}} \int_{{2^k R} \leq y < {2^{k+1} R}} y^{t (\beta - d_\alpha)+ d_\alpha -1} dy \nonumber
  \\
  & = C_\alpha \sum_{k \in \Z} \frac{{(2^k R)}^ {t(\beta - d_\alpha)+ d_\alpha}} {(1 + 2^k R)^{\gamma t}} .\nonumber
\end{align}
Therefore, we have 
\begin{align}
    ||K_{\beta, \gamma}^\alpha||_{L^t(\R, d\mu_\alpha)} \sim \left(\sum_{k \in \Z} \frac{(2^kR)^{(\beta - d_\alpha)t + d_\alpha}}{(1 + 2^kR)^{\gamma t}}\right)^{\frac{1}{t}},\nonumber 
\end{align}
i.e. there exist two constants $B_\alpha$, $C_\alpha$ such that
 \begin{align}
       C_\alpha \left(\sum_{k \in \Z} \frac{(2^kR)^{(\beta - d_\alpha)t + d_\alpha}}{(1 + 2^kR)^{\gamma t}}\right)^{\frac{1}{t}} &\leq ||K_{\beta, \gamma}^\alpha||_{L^t(\R, d\mu_\alpha)} \nonumber
       \\
       &\leq B_\alpha\left(\sum_{k \in \Z} \frac{(2^kR)^{(\beta - d_\alpha)t + d_\alpha}}{(1 + 2^kR)^{\gamma t}}\right)^{\frac{1}{t}}, \nonumber
        \end{align} 
which completes the proof.  \qed

Now using Lemma \ref{lem: 4.1}, we prove the following Theorem.
\begin{theorem}\label{eq:3}
    Let $ 0 < \beta < d_\alpha$ and $\gamma > 0$. If $\phi(r) \leq C r^\nu$ for every $r > 0$, with $\nu < -\beta$, then for $1 < p < \infty$ and $t \in \left( \frac{d_\alpha}{d_\alpha + \gamma - \beta}, \frac{d_\alpha}{d_\alpha - \beta}\right)$, the Dunkl-type Bessel-Riesz operator $I_{\beta, \gamma}^\alpha $ is bounded from  $L^{p, \phi}(\R, d\mu_\alpha)$ to  $L^{q, \psi}(\R, d\mu_\alpha)$ i.e.
    \begin{align*}
        ||I_{\beta, \gamma}^\alpha f||_{L^{q, \psi}(\R, d\mu_\alpha)} \leq C ||K_{\beta, \gamma}^\alpha||_{L^t(\R, d\mu_\alpha)} ||f||_{L^{p, \phi}(\R, d\mu_\alpha)},
    \end{align*}
    where $q = \frac{\nu t' p}{\nu t' + d_\alpha}$ and $\psi(r) = (\phi(r))^{\frac{p}{q}}$.
\end{theorem}
\textbf{Proof.} Let $R > 0$ be given. Then for $f \in L^{p, \phi}(\R, d\mu_\alpha)$, we write
    $$
   I_{\beta, \gamma}^\alpha f (x) = \left( \int_{B(0, R)} + \int_{B^c(0, R)}\right) f(y) \tau_x^\alpha K_{\beta, \gamma}^\alpha(y) d\mu_\alpha(y).
   $$
   This implies
   \begin{align}
      | I_{\beta, \gamma}^\alpha f (x)| &\leq \left( \int_{B(0, R)} + \int_{B^c(0, R)}\right) |f|(y) \tau_x^\alpha K_{\beta, \gamma}^\alpha(y) d\mu_\alpha(y)  \nonumber
      \\
     &= \left( \int_{B(0, R)} + \int_{B^c(0, R)}\right)  \tau_{-x}^\alpha|f|(y) K_{\beta, \gamma}^\alpha(y) d\mu_\alpha(y) = F_1(x) + F_2(x), \label{eq: 4.2c}
   \end{align}
where \begin{align}
F_1(x) := \int_{B(0, R)} \tau_{-x}^\alpha |f|(y) K_{\beta, \gamma}^\alpha(y) d\mu_\alpha(y) \nonumber
    \end{align}
    and
    \begin{align}
F_2(x) := \int_{B^c(0, R)} \tau_{-x}^\alpha |f|(y) K_{\beta, \gamma}^\alpha(y) d\mu_\alpha(y).\nonumber
\end{align}
For the first term, from \eqref{eq: 3.4a}, we have
\begin{align}
    F_1(x) & = \int_{B(0, R)}  \tau_{-x}^\alpha|f|(y) K_{\beta, \gamma}^\alpha(y) d\mu_\alpha(y) \nonumber
    \\
    & = \sum_{k= -\infty}^{-1} \int_{2^k R \leq |y| < 2^{k+1}R} \tau_{-x}^\alpha|f|(y)  \frac{|y|^{\beta - d_\alpha}}{(1 + |y|)^\gamma} d\mu_\alpha(y) \nonumber
    \\
    & \leq \sum_{k= -\infty}^{-1}   \frac{(2^k R)^{\beta - d_\alpha}}{(1 + 2^k R)^\gamma} \int_{2^k R \leq |y| < 2^{k+1}R} \tau_{-x}^\alpha|f|(y)  d\mu_\alpha(y) \nonumber
    \\
    & \leq C M^\alpha f(-x) \sum_{k= -\infty}^{-1} \frac{(2^k R)^{\beta - d_\alpha + \frac{d_\alpha}{t}+ \frac{d_\alpha}{t'}}}{(1 + 2^k R)^\gamma}. \nonumber
    \\
    &= C M^\alpha f(-x) \sum_{k= -\infty}^{-1} \frac{(2^k R)^{\frac{(\beta - d_\alpha)t + d_\alpha}{t}}}{(1 + 2^k R)^\gamma} (2^k R)^ {d_\alpha/t'} \nonumber
\end{align}
By H\"{o}lder inequality for $\frac{1}{t} + \frac{1}{t'} = 1$, we get
\begin{align}
|F_1(x)|  
&\leq C M^\alpha f(-x) \left(\sum_{k= -\infty}^{-1} \frac{(2^k R)^{(\beta - d_\alpha)t + d_\alpha}}{(1 + 2^k R)^{\gamma t}} \right)^{\frac{1}{t}} \left(\sum_{k= -\infty}^{-1} (2^k R)^ {d_\alpha}  \right)^{\frac{1}{t'}} . \nonumber
\end{align}
Since by Lemma \ref{lem: 4.1}, 
\begin{align}
\left(\sum_{k= -\infty}^{-1} \frac{(2^k R)^{(\beta - d_\alpha)t + d_\alpha}}{(1 + 2^k R)^{\gamma t}} \right)^{\frac{1}{t}}  &\leq \left(\sum_{k \in \Z} \frac{(2^k R)^{(\beta - d_\alpha)t + d_\alpha}}{(1 + 2^k R)^{\gamma t}} \right)^{\frac{1}{t}} \nonumber \\
&\leq C ||K_{\beta, \gamma}^\alpha||_{L^t(\R, d\mu_\alpha)},\label{eq: 3.2}
\end{align}
we obtain
\begin{align}
    |F_1(x)| \leq C||K_{\beta, \gamma}^\alpha||_{L^t(\R, d\mu_\alpha)}  M^\alpha f(-x) R^{\frac{d_\alpha}{t'}}. \label{eq: 3.3}
   \end{align}
For the second term $F_2$, we get
\begin{align}
     F_2(x) & = \int_{B^c(0, R)}  \tau_{-x}^\alpha|f|(y) K_{\beta, \gamma}^\alpha(y) d\mu_\alpha(y) \nonumber
    \\
    & = \sum_{k = 0}^{\infty} \int_{2^k R \leq |y| < 2^{k+1}R}  \tau_{-x}^\alpha|f|(y) \frac{|y|^{\beta - d_\alpha}}{(1 + |y|)^\gamma} d\mu_\alpha(y) \nonumber
    \\
    & \leq \sum_{k= 0}^{\infty} \frac{(2^k R)^{\beta - d_\alpha}}{(1 + 2^k R)^\gamma} \int_{2^k R \leq |y| < 2^{k+1}R} \tau_{-x}^\alpha |f|(y)  d\mu_\alpha(y). \nonumber
\end{align}
For each $k = 0, 1, 2, ...,$ we will first estimate 
\begin{align}
    \int_{2^k R \leq |y| < 2^{k+1} R} \tau_{-x}^\alpha |f|(y) d\mu_\alpha(y). \nonumber
\end{align}
Now
\begin{align}
    \int_{2^k R \leq |y| < 2^{k+1} R} \tau_{-x}^\alpha |f|(y) d\mu_\alpha(y) &=  \int_{\mathbb{R}} |f|(y) \tau_{x}^\alpha \chi_{\{2^k R \leq |y| < 2^{k+1} R\}}(y)  d\mu_\alpha(y)  \nonumber
    \\
    & \leq \int_{\mathbb{R}} |f|(y) \tau_{x}^\alpha \chi_{B(0, 2^{k+1}R)}(y) d\mu_\alpha(y) \nonumber
    \\
    &= \int_{\mathbb{R}} |f|(y) (\tau_{x}^\alpha \chi_{B(0, 2^{k+1}R)}(y))^{\frac{1}{p} + \frac{1}{q}} d\mu_\alpha(y), \nonumber
\end{align}
with $p \in (1,\infty)$ and its conjugate $q$ defined through $\frac{1}{p} + \frac{1}{q} = 1$. Then, by applying Hölder’s inequality together with Theorem \ref{tm: 2.1}, we obtain
\begin{align}
  &\int_{2^k R \leq |y| < 2^{k+1} R} \tau_{x}^\alpha |f|(y) d\mu_\alpha(y) \nonumber
  \\
  &\leq  \left(\int_\R |f(y)|^p \tau_{x}^\alpha \chi_{B(0,2^{k+1}R)}(y) d\mu_\alpha(y) \right)^{1/p} \nonumber
    \\
    & \hspace{2.25cm}\times\left(\int_\R \tau_{x}^\alpha \chi_{B(0,2^{k+1}R)}(y) d\mu_\alpha(y) \right)^{1/q}  \nonumber
    \\
    &\leq  \left(\int_\R |f(y)|^p \tau_{x}^\alpha \chi_{B(0,2^{k+1}R)}(y) d\mu_\alpha(y) \right)^{1/p} \nonumber
    \\
    & \hspace{2.25cm}\times\left(\int_\R \chi_{B(0,2^{k+1}R)}(y) d\mu_\alpha(y) \right)^{1/q}  \nonumber
    \\
    & \leq C (2^{k+1} R)^{d_\alpha} \left( \frac{1}{(2^{k+1}R)^{d_\alpha}}\int_\R |f(y)|^p \tau_{x}^\alpha \chi_{B(0,2^{k+1}R)}(y) d\mu_\alpha(y) \right)^{1/p}  \nonumber
    \\
    & \leq C (2^k R)^{d_\alpha} \left( \frac{1}{(2^{k+1}R)^{d_\alpha}}\int_{B(0, 2^{k+1}R)} \tau_{-x}^\alpha |f|^p(y)  d\mu_\alpha(y) \right)^{1/p}.  \nonumber
\end{align}
Thus from the definition of the generalized Dunkl-type Morrey space \ref{eq: 3.2a}, it is clear that
\begin{align}
     \int_{2^k R \leq |y| < 2^{k+1} R } \tau_{-x}^\alpha |f|(y) d\mu_\alpha(y) &  \leq C (2^k R)^{d_\alpha} \phi(2^{k+1} R)||f||_{L^{p,\phi}(\mathbb{R}, d\mu_\alpha)} \nonumber
     \\
     &\leq C \phi(2^k R) (2^k R)^{d_\alpha} ||f||_{L^{p,\phi}(\mathbb{R}, d\mu_\alpha)}. \label{eq: 4.5b}
\end{align}
Using the above inequality \eqref{eq: 4.5b}, we get
\begin{align}
    F_2(x) \leq C ||f||_{L^{p, \phi}(\R, d\mu_\alpha)}\sum_{k= 0}^{\infty} \frac{(2^k R)^{\beta - d_\alpha }}{(1 + 2^k R)^\gamma} (2^k R)^{d_\alpha}\phi(2^k R) .\label{eq: 4.6b}
\end{align}
Using $\phi(r) \leq C r^{\nu}$ and applying H\"{o}lder inequality again, we arrive at
\begin{align}
    F_2(x) & \leq C ||f||_{L^{p, \phi}(\R, d\mu_\alpha)}\sum_{k= 0}^{\infty} \frac{(2^k R)^{\beta - d_\alpha + \frac{d_\alpha}{t}}}{(1 + 2^k R)^\gamma} (2^k R)^{\nu + \frac{d_\alpha}{t'}} \nonumber 
     \\
     & \leq  C ||f||_{L^{p, \phi}(\R, d\mu_\alpha)}\left(\sum_{k= 0}^{\infty} \frac{(2^k R)^{(\beta - d_\alpha)t + d_\alpha}}{(1 + 2^k R)^{\gamma t}}\right)^{\frac{1}{t}} \left(\sum_{k = 0}^{\infty}(2^k R)^{\nu t' + d_\alpha} \right)^{\frac{1}{t'}}. \nonumber 
\end{align}
 Since $\nu t' + d_\alpha < 0$,  Lemma \ref{lem: 4.1} leads to
  \begin{align}
      F_2(x) \leq C ||K_{\beta, \gamma}^\alpha||_{L^{t}(\R, d\mu_\alpha)} ||f||_{L^{p, \phi}(\R, d\mu_\alpha)} R^{\frac{d_\alpha}{t'} + \nu}.\label{eq: 3.6}
  \end{align}
Adding  equations \eqref{eq: 3.3} and \eqref{eq: 3.6}, from \eqref{eq: 4.2c} we get 
\begin{align}
    |I_{\beta, \gamma}^\alpha f(x)| \leq C ||K_{\beta, \gamma}^\alpha||_{L^t(\R, d\mu_\alpha)} \left( M^\alpha f(-x) R^{\frac{d_\alpha}{t'}} + ||f||_{L^{p, \phi}(\R, d\mu_\alpha)} R^{\frac{d_\alpha}{t'} + \nu}\right). \label{eq: 4.8b}
\end{align}
Note that the above inequality \eqref{eq: 4.8b} holds $\forall ~ R > 0$. Suppose that $f \neq 0$ and $M^\alpha f$ is finite everywhere. For each $x \in \R$ take $R_x > 0$ such that $R_x^\nu = \frac{M^\alpha f(-x)}{||f||_{L^{p, \phi}(\R, d\mu_\alpha)}}$. Then by taking $R = R_x$ in \eqref{eq: 4.8b}, we arrive at
\begin{align}
    |I_{\beta, \gamma}^\alpha f(x)| \leq C ||K_{\beta, \gamma}^\alpha||_{L^t(\R, d\mu_\alpha)} ||f||_{L^{p, \phi}(\R, d\mu_\alpha)}^{-\frac{d_\alpha}{\nu t'}}(M^\alpha f(-x))^{1 +
    \frac{d_\alpha}{\nu t'}}. \nonumber
\end{align}
Now for any $ x_o \in \R$ and for any $r >0$, consider 
\begin{align}
 &\left( \int_{B(0,r)}\tau_{x_o}^\alpha|I_{\beta, \gamma}^\alpha f|^q(x) d\mu_\alpha(x) \right)^{\frac{1}{q}} \nonumber\\
 &=  \left(\int_{\R} |I_{\beta, \gamma}^\alpha f|^q(x) \tau_{-x_o}^\alpha \chi_{B(0,r)}(x) d\mu_\alpha(x) \right)^{\frac{1}{q}} \nonumber
    \\ 
 &\leq C ||K_{\beta, \gamma}^\alpha||_{L^t(\R, d\mu_\alpha)} ||f||_{L^{p, \phi}(\R, d\mu_\alpha)}^{ 1 - \frac{p}{q}} \nonumber
  \\
  &\hspace{3cm} \times \left( \int_{\mathbb{R}} |M^\alpha f|^p(-x)  \tau_{-x_o}^\alpha \chi_{B(0,r)}(x) d\mu_\alpha(x) \right)^{\frac{1}{q}}  \nonumber
  \\
   &\leq C ||K_{\beta, \gamma}^\alpha||_{L^t(\R, d\mu_\alpha)} ||f||_{L^{p, \phi}(\R, d\mu_\alpha)}^{ 1 - \frac{p}{q}} \nonumber
  \\
  &\hspace{3cm} \times \left( \int_{B(0,r)} \tau_{x_o}^\alpha |M^\alpha f|^p(-x) d\mu_\alpha(x) \right)^{\frac{1}{q}} , \nonumber
\end{align}
 where $ q =\frac{\nu t' p}{\nu t' + d_\alpha}$. Dividing both sides by $(\phi(r))^{\frac{p}{q}} r^{\frac{d_\alpha}{q}}$, we obtain
 \begin{align}
\frac{\left(\int_{B(0,r)}\tau_{x_o}^\alpha |I_{\beta,\gamma}^\alpha f|^q(x)d\mu_\alpha(x)\right)^{\frac{1}{q}}}{\psi(r)r^{\frac{d_\alpha}{q}} }
 &\leq C ||K_{\beta, \gamma}^\alpha||_{L^t(\R, d\mu_\alpha)}
||f||_{L^{p, \phi}(\R, d\mu_\alpha)}^{1 - \frac{p}{q}} \nonumber
\\
& \hspace{1.25cm}\times \frac{\left(\int_{B(0, r)}\tau_{x_o}^\alpha |M^\alpha f|^p(-x) d\mu_\alpha(x)\right)^{\frac{1}{q}}}{(\phi(r))^{\frac{p}{q}}r^{\frac{d\alpha}{q}}}, \nonumber
\end{align}
where $\psi(r)= (\phi(r))^{\frac{p}{q}}$. Now, we take the supremum over $r > 0$ and $x_0 \in \mathbb{R}$ to have
\begin{align}
        ||I_{\beta, \gamma }^\alpha f||_{L^{q, \psi}(\R, d\mu_\alpha)} \leq C ||K_{\beta, \gamma}^\alpha||_{L^t(\R, d\mu_\alpha)} ||f||_{L^{p, \phi}(\R, d\mu_\alpha)}^{1 - \frac{p}{q}} ||M^\alpha f ||_{L^{p, \phi}(\R, d\mu_\alpha)}^{\frac{p}{q}}, \nonumber
\end{align}
which implies 
\begin{align}
     ||I_{\beta, \gamma }^\alpha||_{L^{q, \psi}(\R, d\mu_\alpha)} \leq C ||K_{\beta, \gamma}^\alpha||_{L^t(\R, d\mu_\alpha)} ||f||_{L^{p, \phi}(\R, d\mu_\alpha)}.   \nonumber   
 \end{align} 
 This completes the proof.\qed
 
 Now by using Lemma \ref{lem: 2.4}, we wish to obtain a more general result for the boundedness of $I_{\beta, \gamma}^\alpha$.
 \begin{theorem} \label{eq:4}
      Let $ 0 < \beta < d_\alpha$ and $\gamma > 0$. If $\phi(r) \leq C r^\nu$ for every $r > 0$, with $\nu < -\beta$, $1 < p < \infty$ and $ \frac{d_\alpha}{d_\alpha + \gamma - \beta} < s \leq t < \frac{d_\alpha}{d_\alpha - \beta}$ and $1 \leq s \leq t$,  then for all $f \in L^{p, \phi}(\R, d\mu_\alpha)$ we get
    \begin{align}
        ||I_{\beta, \gamma}^\alpha f||_{L^{q, \psi}(\R, d\mu_\alpha)} \leq C ||K_{\beta, \gamma}^\alpha||_{L^{s, t}(\R, d\mu_\alpha)} ||f||_{L^{p, \phi}(\R, d\mu_\alpha)},
    \end{align}
    where $q = \frac{\nu t' p}{\nu t' + d_\alpha}$ and $\psi(r) = (\phi(r))^{\frac{p}{q}}$. 
 \end{theorem}
 \textbf{Proof.} As in the proof of Theorem \ref{eq:3}, for given $R > 0$ we write
   $$ 
   |I_{\beta, \gamma}^\alpha f (x)| \leq F_1(x) + F_2(x) ,
   $$
where \begin{align}
    F_1(x) := \int_{B(0, R)} \tau_{-x}^\alpha |f|(y)  K_{\beta, \gamma}^\alpha(y) d\mu_\alpha(y) \nonumber
    \end{align}
    and
    \begin{align}
    F_2(x) := \int_{B^c(0, R)} \tau_{-x}^\alpha |f|(y) K_{\beta, \gamma}^\alpha(y) d\mu_\alpha(y). \nonumber
\end{align}
As in Theorem \ref{eq:3}, using dyadic decomposition and  by H\"{o}lder inequality for $\frac{1}{s} + \frac{1}{s'} =1$, we get
\begin{align}
F_1(x) &\leq C M^\alpha f(-x) \left(\sum_{k= -\infty}^{-1} \frac{(2^k R)^{(\beta - d_\alpha)s + d_\alpha}}{(1 + 2^k R)^{\gamma s}} \right)^{\frac{1}{s}} \left(\sum_{k= -\infty}^{-1} (2^k R)^ {d_\alpha}  \right)^{\frac{1}{s'}} . \nonumber
\end{align}
Replacing $t$ by $s$ in \eqref{eq: 3.2} where $1 \leq s \leq t$, we get
\begin{align}
    F_1(x) &\leq C  M^\alpha f(-x) \left( \int_{B(0, R)} (K_{\beta, \gamma}^\alpha (y))^s d\mu_\alpha(y) \right)^{\frac{1}{s}} R^{\frac{d_\alpha}{s'}}  \nonumber\\ 
   & \leq C M^\alpha f(-x) ||K_{\beta, \gamma}^\alpha||_{L^{s, t}(\R, d\mu_\alpha)} R^{\frac{d_\alpha}{t'}}, \label{eq: 3.8}
   \end{align}
   since 
   \begin{align}
       \int_{2^k R \leq |y| < 2^{k+1 R}} (K_{\beta, \gamma})^s(y) d\mu_\alpha(y) \sim \frac{(2^k R)^{s(\beta - d_\alpha) +d_\alpha}}{(1 + 2^k R)^{\gamma s}}. \nonumber
   \end{align}
Now, for the second term $F_2$, proceeding similarly as in the proof of Theorem \ref{eq:3}, we get from \eqref{eq: 4.6b}
\begin{align}
    F_2(x)  &\leq C ||f||_{L^{p, \phi}(\mathbb{R}, d\mu_\alpha)}  \sum_{k= 0}^{\infty} \frac{(2^k R)^{\beta - d_\alpha}}{(1 + 2^k R)^\gamma} \phi(2^kR) (2^k R)^{d_\alpha} \nonumber
     \\
     & \leq C ||f||_{L^{p, \phi}(\R, d\mu_\alpha)}\sum_{k= 0}^{\infty} \frac{(2^k R)^{\beta - d_\alpha} \phi(2^k R)(2^k R)^{d_\alpha}}{(1 + 2^k R)^\gamma} \nonumber
     \\ 
     & \hspace{4.5cm}\times\frac{\left(\int_{2^k R \leq |y| < 2^{k + 1} R} d\mu_\alpha(y)\right)^{\frac{1}{s}}}{(2^k R)^{\frac{d_\alpha}{s}}} , \label{eq: 4.12}
     \end{align}
     where we have used the fact that  $\left( \int_{2^k R \leq |y| < 2^{k+1} R} d\mu_\alpha(y) \right)^{\frac{1}{s}} \sim (2^k R)^{\frac{d_\alpha}{s}}$. Now following similarly as in the proof of Lemma \ref{lem: 4.1}, we can show that
     \begin{align}
         &\left(\int_{2^k R \leq |y| < 2^{k+1}R}(K_{\beta, \gamma}^\alpha)^s(y) d\mu_\alpha(y)\right)^\frac{1}{s} \nonumber
         \\
         &\sim \frac{(2^k R)^{\beta - d_\alpha +\frac{d_\alpha}{s}}}{(1 + 2^k R)^\gamma} \nonumber
         \\
         &\geq C \frac{(2^k R)^{\beta - d_\alpha}}{(1 + 2^k R)^\gamma} \left( \int_{2^k R \leq |y| < 2^{k + 1}R} d\mu_\alpha(y) \right)^\frac{1}{s} . \label{eq: 24}
         \end{align}
 Thus substituting \eqref{eq: 24} in \eqref{eq: 4.12}, we obtain
     \begin{align}
     F_2(x) &\leq C ||f||_{L^{p, \phi}(\R, d\mu_\alpha)}\sum_{k= 0}^{\infty} {(2^k R)^\frac{d_\alpha}{t'} \phi(2^k R)} \nonumber
     \\
   &\hspace{4.5cm} \times \frac{ \left(\int_{2^k R \leq |y| < 2^{k+1}R}(K_{\beta, \gamma}^\alpha)^s(y) d\mu_\alpha(y)\right)^\frac{1}{s}}{(2^k R)^{\frac{d_\alpha}{s} -\frac{d_\alpha}{t}}} \nonumber
     \\
     &\leq C ||K_{\beta, \gamma}^\alpha||_{L^{s,t}(\R, d\mu_\alpha)} ||f||_{L^{p, \phi}(\R, d\mu_\alpha)} \sum_{k=0}^{\infty} (2^k R)^{\nu + \frac{d_\alpha}{t'}}.\nonumber
  \end{align}
  Using the assumption $\phi(r) \leq C r^\nu$ and  $\nu + \frac{d_\alpha}{t'} < 0$ , we have 
\begin{align}
  F_2(x) & \leq C ||K_{\beta, \gamma}^\alpha||_{L^{s,t}(\R, d\mu_\alpha)} ||f||_{L^{p, \phi}(\R, d\mu_\alpha)} R^{ \frac{d_\alpha}{t'} + \nu}.\label{eq: 3.10}
  \end{align}
Now adding  equations \eqref{eq: 3.8} and \eqref{eq: 3.10}, we get 
\begin{align}
    |I_{\beta, \gamma}^\alpha f(x)| \leq C ||K_{\beta, \gamma}^\alpha||_{L^{s,t}(\R, d\mu_\alpha)}\left( M^\alpha f(-x) R^{\frac{d_\alpha}{t'}} + ||f||_{L^{p, \phi}(\R, d\mu_\alpha)} R^{\frac{d_\alpha}{t'} + \nu}\right). \nonumber
\end{align}
Now proceeding similarly as in the proof of Theorem \ref{eq:3}, we arrive at
\begin{align}
    ||I_{\beta, \gamma }^\alpha f||_{L^{q, \psi}(\R, d\mu_\alpha)} \leq C ||K_{\beta, \gamma}^\alpha||_{L^{s,t}(\R, d\mu_\alpha)} ||f||_{L^{p, \phi}(\R, d\mu_\alpha)}, \nonumber
  \end{align}
  which completes the proof.       \qed 
  
 By Lemma \ref{lem: 2.4}, we immediately get that Theorem \ref{eq:4} implies Theorem \ref{eq:3}:
 \begin{align}
    ||I_{\beta, \gamma }^\alpha f||_{L^{q, \psi}(\R, d\mu_\alpha)} & \leq C ||K_{\beta, \gamma}^\alpha||_{L^{s,t}(\R, d\mu_\alpha)} ||f||_{L^{p, \phi}(\R, d\mu_\alpha)} \nonumber
    \\
    &\leq C ||K_{\beta, \gamma}^\alpha||_{L^{t}(\R, d\mu_\alpha)} ||f||_{L^{p, \phi}(\R, d\mu_\alpha)} . \nonumber
 \end{align}
 
 We still wish to obtain a better estimate. First, we prove the following lemma, which presents that the kernel $K_{\beta, \gamma}^\alpha$ belongs to the generalized Morrey space $L^{s, \omega}(\mathbb{R}, d\mu_\alpha)$ for some $s \geq 1$ and some function $\omega$.
 \begin{lemma}
     If s $\geq 1 , \gamma > 0$ and $ \omega : \R^{+} \to \R^{+}$ with $ \omega(r) \geq C r^{\beta - d_\alpha}$ for every $r > 0$, then $K_{\beta, \gamma}^\alpha \in L^{s, \omega}(\R, d\mu_\alpha)$, where $d_\alpha - \frac{d_\alpha}{s} < \beta < d_\alpha$.
 \end{lemma}

 \textbf{Proof.} In order to prove this lemma, we will derive the following:
 \begin{align}
     \left( \int_{B(0, R)} \tau_x^\alpha (K_{\beta, \gamma }^\alpha)^s(y) d\mu_\alpha(y) \right)^{\frac{1}{s}}
      &\leq C R^{\beta - d_\alpha + \frac{d_\alpha}{s}} . \label{eq: 281}
 \end{align}
 The above inequality was proved in \cite{sps} for the kernel $K_\beta^\alpha$. But for the sake of completeness, we provide the proof here.
  For $x \in \mathbb{R}$, consider
 \begin{align}
     \left( \int_{B(0, R)} \tau_x^\alpha (K_{\beta, \gamma }^\alpha)^s(y) d\mu_\alpha(y) \right)^{\frac{1}{s}}. \nonumber
 \end{align}
 First, we take $|x| > 2 R$. Then we have $|x \pm y| > R$, which implies $d(x, y) \geq R$ where $d(x, y)= min \{ |x-y| ,|x + y|\}$.
 We shall show that 
 \begin{align}
\tau_x^\alpha K_{\beta, \gamma}^\alpha(y) \leq d(x, y)^{\beta - d_\alpha}. \nonumber
\end{align}
Since, by \eqref{eq: 2.2}
 \begin{align}
     \tau_x^\alpha K_{\beta, \gamma}^ \alpha (y) & = \int_\R \frac{(\sqrt{x^2 + y^2 - 2x \eta })^{\beta - d_\alpha}}{(1 + \sqrt{{x^2 + y^2 - 2x \eta }})^\gamma} d\mu_y^\alpha(\eta)  \nonumber
     \\
     & \leq \int_\R (\sqrt{x^2 + y^2 - 2x\eta})^{\beta- d_\alpha} d\mu_y^\alpha(\eta), \nonumber
 \end{align}
and Amri and Sifi \cite{amri2012riesz} proved that 
 \begin{align}
     min\{ |x - y| , |x + y| \} \leq A(x,y, \eta) \leq max\{|x - y| , |x+ y| \} ,
     \nonumber
 \end{align}
 $ \forall~ x, y \in \R , \forall ~\eta \in [-y, y]$, where $A(x, y, \eta) = \sqrt{x^2 + y^2 - 2x \eta}$, we have 
     \begin{align}
         \tau_x^\alpha K_{\beta, \gamma}^\alpha(y) \leq \int_\R d(x, y)^{\beta - d_\alpha} d\mu_y^\alpha(\eta) = d(x,y)^{\beta - d_\alpha}. \nonumber
         \end{align}
      Using above estimates we immediately get $\tau_x^\alpha(K_{\beta, \gamma}^\alpha(y))^s \leq d(x,y)^{(\beta - d_\alpha)s}$.\\
   We have $d(x, y) \geq R$ because $|x \pm y| > R$ for $|x| > 2R$ and $|y| < R$, this implies 
   \begin{align}
       \left(\int_{B(0, R)} \tau_x^\alpha (K_{\beta, \gamma}^\alpha)^s(y) d\mu_\alpha(y) \right)^{\frac{1}{s}} &\leq \left( \int_{B(0,R)} d(x,y)^{(\beta -d_\alpha)s} d\mu_\alpha(y) \right)^{\frac{1}{s}} \nonumber \\
       &\leq C R^{\beta - d_\alpha + \frac{d_\alpha}{s}}. \label{eq: 27}
   \end{align}
   Next we consider the case $|x| \leq 2R$. We have $\tau_x^\alpha K_{\beta, \gamma}^\alpha(y) \geq 0$ as $K_{\beta, \gamma}^\alpha$ is positive radial function. Then
   \begin{align}
       &\left( \int_{B(0, R)} \tau_x^\alpha(K_{\beta, \gamma}^\alpha(y))^s d\mu_\alpha(y) \right)^\frac{1}{s}\nonumber
       \\
       \leq \left( \int_{B(0, 2 R)} \tau_x^\alpha(K_{\beta, \gamma}^\alpha(y))^s d\mu_\alpha(y) \right)^\frac{1}{s} \nonumber 
       \\
        &= \left( \int_{\mathbb{R}} (K_{\beta, \gamma}^\alpha(y))^s \tau_{-x}^\alpha \chi_{B(0, 2 R)}(y) d\mu_\alpha(y) \right)^\frac{1}{s} \nonumber
        \\
       & = \left( \int_{B(-x, 2 R)} (K_{\beta, \gamma}^\alpha(y))^s \tau_{-x}^\alpha \chi_{B(0, 2 R)}(y) d\mu_\alpha(y) \right)^\frac{1}{s}. \nonumber
    \end{align}
    By using Lemma \ref{lem: 2.2} and the fact that $B(-x, 2 R) = B(0, |x| + 2R) \subseteq B(0, 4R)$ in the above inequality, we get 
    \begin{align}
        \left( \int_{B(0, R)} \tau_x^\alpha(K_{\beta, \gamma}^\alpha(y))^s d\mu_\alpha(y) \right)^\frac{1}{s} &\leq \left( \int_{B(0, |x| + 2 R)} (K_{\beta, \gamma}^\alpha(y))^s d\mu_\alpha(y) \right)^\frac{1}{s} \nonumber
        \\
        &\leq \left( \int_{B(0, 4R)}\frac{|y|^{(\beta - d_\alpha)s}}{(1 + |y|)^{\gamma s}} d\mu_\alpha (y) \right)^\frac{1}{s}\nonumber
        \\
        &\leq C R^{\beta - d_\alpha + \frac{d_\alpha}{s}}  \label{eq: 28},
    \end{align}
    since $s < \frac{d_\alpha}{d_\alpha - \beta}$. Thus combining \eqref{eq: 27} and \eqref{eq: 28}, we obtain \eqref{eq: 281}. \\
    Hence, we have for all $x \in \R$ and $R > 0$ (using $\omega(R) \geq C R^{\beta - d_\alpha}$)
    \begin{align}
        \left( \int_{B(0, R)} \tau_x^\alpha (K_{\beta, \gamma}^\alpha(y))^s d\mu_\alpha(y) \right)^\frac{1}{s} \leq C \omega(R) R^\frac{d_\alpha}{s}. \nonumber
    \end{align}
    By dividing both sides of the above inequality by $\omega(R) R^\frac{d_\alpha}{s}$, we get
    \begin{align}
        \frac{1}{\omega(R) R^\frac{d_\alpha}{s}} \left( \int_{B(0, R)} \tau_x^\alpha (K_{\beta, \gamma}^\alpha(y))^s d\mu_\alpha(y) \right)^\frac{1}{s} \leq C. \nonumber
    \end{align}
   Now, we take supremum over $R > 0$ and $x \in \R$ to get 
    \begin{align}
   ||K_{\beta, \gamma}^\alpha||_{L^{s, \omega}(\R, d\mu_\alpha)} =\sup_{\substack{{R>0}\\{x \in \R}}}  \frac{1}{\omega(R) R^\frac{d_\alpha}{s}} \left( \int_{B(0, R)} \tau_x^\alpha (K_{\beta, \gamma}^\alpha(y))^s d\mu_\alpha(y) \right)^\frac{1}{s} < \infty, \label{eq: 29} 
    \end{align}
  which implies $K_{\beta, \gamma}^\alpha \in L^{s, \omega}(\R, d\mu_\alpha)$. \qed  
  \begin{theorem}
 Assume that $\omega : \R^{+} \to \R^{+}$ satisfies the doubling condition and $ C r^{\beta - d_\alpha} \leq \omega(r) \leq C r^{-\beta}$ for every $r > 0$. Let $\gamma > 0$ and $ 0 < \beta < d_\alpha$. Further assume that $\frac{d_\alpha}{d_\alpha + \gamma -\beta} < s < \frac{d_\alpha}{d_\alpha - \beta}$, $s \geq 1$. If $\phi(r) \leq C r^\nu $ for every $r > 0$, where $\nu < - \beta < -d_\alpha - \nu$, then for all $f \in L^{p, \phi}(\R, d\mu_\alpha)$, we have
      \begin{align}
          ||I_{\beta, \gamma}^\alpha f||_{L^{q, \psi}(\R, d\mu_\alpha)} \leq C ||K_{\beta, \gamma}^\alpha||_{L^{s, \omega}(\R, d\mu_\alpha)} ||f||_{L^{p, \phi}(\R, d\mu_\alpha)}, \nonumber
    \end{align}
      with $1 < p < \infty, q = \frac{\nu p}{\nu + d_\alpha -\beta }$ and $\psi(r) = (\phi(r))^\frac{p}{q}$. 
  \end{theorem}
  \textbf{Proof.} Let $R > 0$ be given. As in the proof of Theorem \ref{eq:3}, for $f \in L^{p, \phi}(\R, d\mu_\alpha)$ we write
   $$ 
   |I_{\beta, \gamma}^\alpha f (x)| \leq  F_1(x) + F_2(x) ,
   $$
where \begin{align}
    F_1(x) := \int_{B(0, R)} \tau_{-x}^\alpha  |f|(y)  K_{\beta, \gamma}^\alpha(y) d\mu_\alpha(y) \nonumber
    \end{align}
    and
    \begin{align}
    F_2(x) := \int_{B^c(0, R)}  \tau_{-x}^\alpha |f|(y)  K_{\beta, \gamma}^\alpha(y) d\mu_\alpha(y). \nonumber
\end{align}
 For the first term $F_1$, proceeding similarly as in the proof of Theorem \ref{eq:4}, we get
\begin{align}
    F_1(x) &\leq C  M^\alpha f(-x) \left( \int_{B(0, R)} (K_{\beta, \gamma}^\alpha)^s (y) d\mu_\alpha(y) \right)^{\frac{1}{s}} R^{\frac{d_\alpha}{s'}}  \nonumber
    \\ 
    &\leq C  M^\alpha f(-x) \omega(R) \frac{1}{\omega(R)}\left(\frac{1}{R^{d_\alpha}} \int_{B(0, R)} (K_{\beta, \gamma}^\alpha)^s (y) d\mu_\alpha(y) \right)^{\frac{1}{s}} R^{d_\alpha} \nonumber
    \\
   & \leq C M^\alpha f(-x) ||K_{\beta, \gamma}^\alpha||_{L^{s, \omega}(\R, d\mu_\alpha)} R^{d_\alpha} \omega(R) \nonumber
    \\ 
     &\leq C M^\alpha f(-x) ||K_{\beta, \gamma}^\alpha||_{L^{s, \omega}(\R, d\mu_\alpha)} R^{d_\alpha -\beta}  \label{eq:d"},
   \end{align}
   since $\omega(R) \leq C R^{- \beta}$.\\
Again, for the second term $F_2$ proceeding similarly as in the proof of Theorem \ref{eq:4}, we get
\begin{align}
    F_2(x)  & \leq C ||f||_{L^{p, \phi}(\R, d\mu_\alpha)}\sum_{k= 0}^{\infty} \frac{(2^k R)^\beta \phi(2^k R)}{(1 + 2^k R)^\gamma} \frac{\left(\int_{2^k R \leq |y| < 2^{k + 1} R} d\mu_\alpha(y)\right)^{\frac{1}{s}}}{(2^k R)^{\frac{d_\alpha}{s}}}. \nonumber
     \end{align}
    By using equation \eqref{eq: 24} we have,
     \begin{align}
     F_2(x) \leq C ||f||_{L^{p, \phi}(\R, d\mu_\alpha)}&\sum_{k= 0}^{\infty} \frac{(2^k R)^\beta \phi(2^k R)} {(2^k R)^{\beta - d_\alpha}} \nonumber
     \\
     & \hspace{1.5cm} \times\frac{ \left(\int_{2^k R \leq |y| < 2^{k+1}R}(K_{\beta,\gamma}^\alpha)^s(y)d\mu_\alpha(y)\right)^\frac{1}{s}}{(2^k R)^{\frac{d_\alpha}{s}}} ,\nonumber
  \end{align}
since $\phi(r) \leq C r^\nu$ and $\omega(r) \leq C r^{- \beta}$ for every $r > 0$, using this we have 
\begin{align}
   F_2(x) \leq C ||f||_{L^{p, \phi}(\R, d\mu_\alpha)}&\sum_{k= 0}^{\infty} (2^k R)^{d_\alpha + \nu - \beta} \nonumber
   \\
   &\hspace{1.5cm} \times\frac{ \left(\int_{2^k R \leq |y| < 2^{k+1}R}(K_{\beta,\gamma}^\alpha)^s(y)d\mu_\alpha(y)\right)^\frac{1}{s}}{{\omega(2^k R)}{(2^k R)^{\frac{d_\alpha}{s}}}}. \nonumber
\end{align}
Now, using \eqref{eq: 29}, we get
\begin{align}
  F_2(x) &\leq C ||K_{\beta, \gamma}^\alpha||_{L^{s,\omega}(\R, d\mu_\alpha)} ||f||_{L^{p, \phi}(\R, d\mu_\alpha)} \sum_{k=0}^{\infty} (2^k R)^{\nu + d_\alpha -\beta }. \nonumber
\end{align}
 Since $\nu + d_\alpha - \beta < 0$, this implies
  \begin{align}
      F_2(x) \leq C ||K_{\beta, \gamma}^\alpha||_{L^{s,\omega}(\R, d\mu_\alpha)} ||f||_{L^{p, \phi}(\R, d\mu_\alpha)} R^{d_\alpha + \nu - \beta }.\label{eq: e"}
  \end{align}
Adding  equations \eqref{eq:d"} and \eqref{eq: e"}, we get 
\begin{align}
    |I_{\beta, \gamma}^\alpha f(x)| \leq C ||K_{\beta, \gamma}^\alpha||_{L^{s,\omega}(\R, d\mu_\alpha)}\left( M^\alpha f(-x) R^{d_\alpha - \beta} + ||f||_{L^{p, \phi}(\R, d\mu_\alpha)} R^{d_\alpha -\beta + \nu}\right). \label{eq: 4.20}
\end{align}
Note that the above inequality \eqref{eq: 4.20} holds $ \forall~~ R > 0$. Suppose that $f \neq 0$ and $M^\alpha f$ is finite everywhere. For each $x \in \R$ choose $R_x > 0$ such that $R_x^\nu = \frac{M^\alpha f(-x)}{||f||_{L^{p, \phi}(\R, d\mu_\alpha)}}$. Then by taking $R = R_x$ in \eqref{eq: 4.20}, we arrive at
\begin{align}
    |I_{\beta, \gamma}^\alpha f(x)| \leq C ||K_{\beta, \gamma}^\alpha||_{L^{s,\omega}(\R, d\mu_\alpha)} ||f||_{L^{p, \phi}(\R, d\mu_\alpha)}^{\frac{ \beta - d_\alpha}{\nu}}(M^\alpha f(-x))^{1 +\frac{d_\alpha - \beta}{\nu}}. \nonumber
\end{align}
Now for each $ x_o \in \R$ and for each $r > 0$, consider
\begin{align}
 \left( \int_{B(0,r)}\tau_{x_o}^\alpha|I_{\beta, \gamma}^\alpha f|^q(x) d\mu_\alpha(x) \right)^{\frac{1}{q}} &=  \left(\int_{\R} |I_{\beta, \gamma}^\alpha f|^q(x) \tau_{-x_o}^\alpha \chi_{B(0,r)}(x) d\mu_\alpha(x) \right)^{\frac{1}{q}} \nonumber
    \\ 
 &\leq C ||K_{\beta, \gamma}^\alpha||_{L^{s, \omega}(\R, d\mu_\alpha)} ||f||_{L^{p, \phi}(\R, d\mu_\alpha)}^{ 1 - \frac{p}{q}} \nonumber
  \\ & \hspace{0.50 cm}\times \left( \int_{\mathbb{R}} |M^\alpha f|^p(-x) \tau_{-x_o}^\alpha \chi_{B(0,r)}d\mu_\alpha(x) \right)^{\frac{1}{q}}  \nonumber
  \\
  &\leq C ||K_{\beta, \gamma}^\alpha||_{L^{s, \omega}(\R, d\mu_\alpha)} ||f||_{L^{p, \phi}(\R, d\mu_\alpha)}^{ 1 - \frac{p}{q}} \nonumber
  \\ & \hspace{1 cm}\times \left( \int_{B(0,r)} \tau_{x_o}^\alpha |M^\alpha f|^p(-x) d\mu_\alpha(x) \right)^{\frac{1}{q}}, \nonumber
\end{align}
where $ q =\frac{\nu p}{\nu + d_\alpha - \beta }$. Dividing both sides by $(\phi(r))^{\frac{p}{q}} r^{\frac{d_\alpha}{q}}$, we obtain
 \begin{align}
\frac{\left(\int_{B(0,r)}\tau_{x_o}^\alpha |I_{\beta,\gamma}^\alpha f|^q(x)d\mu_\alpha(x)\right)^{\frac{1}{q}}}{\psi(r)r^{\frac{d_\alpha}{q}} }
 &\leq C ||K_{\beta, \gamma}^\alpha||_{L^{s,\omega}(\R, d\mu_\alpha)}
||f||_{L^{p, \phi}(\R, d\mu_\alpha)}^{1 - \frac{p}{q}}  \nonumber
\\ 
   & \hspace{1.25 cm}\times \frac{\left(\int_{B(0, r)}\tau_{x_o}^\alpha |M^\alpha f|^p(-x) d\mu_\alpha(x)\right)^{\frac{1}{q}}}{(\phi(r))^{\frac{p}{q}}r^{\frac{d\alpha}{q}}}, \nonumber
\end{align}
where $\psi(r)= (\phi(r))^{\frac{p}{q}}$. Now, we take the supremum over $r > 0$ and $x_o \in \R$ to have
\begin{align}
        ||I_{\beta, \gamma }^\alpha f||_{L^{q, \psi}(\R, d\mu_\alpha)} \leq C ||K_{\beta, \gamma}^\alpha||_{L^{s, \omega}(\R, d\mu_\alpha)} ||f||_{L^{p, \phi}(\R, d\mu_\alpha)}^{1 - \frac{p}{q}} ||M^\alpha f ||_{L^{p, \phi}(\R, d\mu_\alpha)}^{\frac{p}{q}}, \nonumber
\end{align}
which implies
\begin{align}
    ||I_{\beta, \gamma }^\alpha f||_{L^{q, \psi}(\R, d|\mu_\alpha)} \leq C ||K_{\beta, \gamma}^\alpha||_{L^{s,\omega}(\R, d\mu_\alpha)} ||f||_{L^{p, \phi}(\R, d\mu_\alpha)}. \nonumber
\end{align} 
This completes the proof. \qed


\section{ Inequalities for generalized Dunkl-Type  Bessel-Riesz Operator in Generalized Dunkl-type Morrey spaces}
\label{Sec:5}
In this section, we prove the generalized  Dunkl-type Bessel-Riesz operator $I_{\TR,\gamma}^\alpha$ is bounded on the generalized  Dunkl-type Morrey-space \eqref{eq: 3.2a}.
First, we define the generalized Dunkl-type Bessel-Riesz operator $I_{\TR, \gamma}^\alpha$ by 
\begin{align}
    I_{\TR, \gamma}^\alpha f(x) := \int_\R f(y)  \tau_x^\alpha \frac{\TR(|y|)}{(1+ |y|)^\gamma} d\mu_\alpha(y), \label{eq:34}
\end{align}
where $\gamma \geq 0$ , $\TR : \R^{+} \to \R^{+}$ and $\TR$ satisfies the doubling condition \eqref{eq: 3.3a} along with the condition 
\begin{align}
      \int_0^\infty \frac{\TR(t)}{t^{\gamma - d_\alpha + 1}} dt < \infty.  \label{eq: 35}
\end{align}
Here $K_{\TR, \gamma}^\alpha(x) = \frac{\TR(|x|)}{(1 + |x|)^\gamma}$ is known as the generalized Dunkl-type Bessel-Riesz kernel. If $\TR(t) = t^{\beta - d_\alpha}$, $\gamma < \beta < d_\alpha$, then we have the Dunkl-type Bessel-Riesz operator $I_{\beta, \gamma}^\alpha = I_{\TR, \gamma}^\alpha$ with Dunkl-type Bessel-Riesz kernel $K_{\TR, \gamma}^\alpha = K_{\beta, \gamma}^\alpha$.
\begin{theorem}
    Let $\gamma > 0$. Further assume that $\TR$ and $\phi$ satisfy the doubling condition \eqref{eq: 3.3a}. Let $\phi$ be surjective and for $1 < p < q < \infty$, it satisfies
    \begin{align}
        \phi(r) \int_0^r \frac{\TR(t)}{t^{\gamma - d_\alpha + 1}} dt + \int_r^\infty \frac{\TR(t) \phi(t)}{t^{\gamma - d_\alpha +1}} dt \leq C (\phi(r))^{\frac{p}{q}}, \label{eq:37}
    \end{align}
    for all $r > 0$. Then the generalized Dunkl-type Bessel-Riesz operator $I_{\TR, \gamma}^\alpha$ is bounded from $L^{p, \phi}(\mathbb{R}, d\mu_\alpha)$ to $L^{q, \psi}(\mathbb{R}, d\mu_\alpha)$ i.e.
    \begin{align}
        ||I^\alpha_{\TR, \gamma}f ||_{L^{q, \psi}(\R, d\mu_\alpha)} \leq C ||f||_{L^{p, \phi}(\R, d\mu_\alpha)}, \label{eq: 38}
    \end{align}
    where $\psi(r) = (\phi(r))^{\frac{p}{q}}$.
    \end{theorem}
    \textbf{Proof.} Let $R > 0$ be given. Then for $f \in L^{p, \phi}(\R, d\mu_\alpha)$, we write
   \begin{align}
   |I_{\TR, \gamma}^\alpha f (x)| &\leq \left(\int_{B(0,R)} + \int_{B^c(0,R)}\right) |f|(y) \tau_x^\alpha \frac{\TR(|y|)}{(1+|y|)^\gamma} d\mu_\alpha(y) \nonumber
   \\
   &= \left(\int_{B(0,R)} + \int_{B^c(0,R)}\right) \tau_{-x}^\alpha  |f|(y) \frac{\TR(|y|)}{(1+|y|)^\gamma} d\mu_\alpha(y) \nonumber
   \\
   &= F_{1,\TR}(x) + F_{2, \TR}(x) ,\nonumber
 \end{align}
where
    \begin{align}
    F_{1, \TR}(x) := \int_{B(0, R)} \tau_{-x}^\alpha  |f|(y) \frac{\TR(|y|)}{(1 + |y|)^\gamma} d\mu_\alpha(y) \nonumber
    \end{align}
    and
    \begin{align}
    F_{2,\TR}(x) := \int_{B^c(0, R)} \tau_{-x}^\alpha |f|(y) \frac{\TR(|y|)}{(1 + |y|)^\gamma}  d\mu_\alpha(y). \nonumber
    \end{align}
 For the first term $F_{1, \TR}$, using doubling condition \eqref{eq: 3.3a} for $\TR$ and from \eqref{eq: 3.4a}, we get
\begin{align}
    F_{1, \TR}(x) & = \int_{B(0, R)}  \tau_{-x}^\alpha |f|(y)  \frac{\TR(|y|)}{(1 + |y|)^\gamma}d\mu_\alpha(y) \nonumber
    \\
      & \leq \int_{B(0, R)}  \tau_{-x}^\alpha |f|(y)  \frac{\TR(|y|)}{(|y|)^\gamma}d\mu_\alpha(y) \nonumber
    \\
    & = \sum_{k= -\infty}^{-1} \int_{2^k R \leq |y| < 2^{k+1}R} \tau_{-x}^\alpha |f|(y)  \frac{\TR(|y|)}{( |y|)^\gamma} d\mu_\alpha(y) \nonumber
    \\
    & \leq C\sum_{k= -\infty}^{-1} \frac{\TR(2^k R)}{( 2^k R)^\gamma}\int_{2^k R \leq |y| < 2^{k+1}R} \tau_{-x}^\alpha |f|(y)  d\mu_\alpha(y) \nonumber
    \\
    & \leq C M^\alpha f(-x) \sum_{k= -\infty}^{-1} \frac{\TR(2^k R)}{(2^k R)^{\gamma - d_\alpha}}. \label{eq: 5.5}
\end{align}
Since 
\begin{align}
    \int_{2^k R}^{2^{k+1}R} \frac{\TR(t)}{t^{\gamma - d_\alpha +1}} dt \geq C \frac{\TR(2^k R)}{(2^k R)^{\gamma -d_\alpha +1}} 2^k R \geq C \frac{\TR(2^k R)}{(2^k R)^{\gamma - d_\alpha}}, \nonumber
\end{align}
 using \eqref{eq:37} and the above estimate in \eqref{eq: 5.5}, we have
\begin{align}
    F_{1, \TR}(x) &\leq C M^\alpha f(-x) \sum_{k=-\infty}^{-1} \int_{2^k R}^{2^{k+1}R} \frac{\TR(t)}{t^{\gamma - d_\alpha +1}} dt \nonumber
    \\
    &= C M^\alpha f(-x) \int_0^R \frac{\TR(t)}{t^{\gamma - d_\alpha +1}} dt \nonumber
    \\ 
    & \leq C M^\alpha f(-x) (\phi(R))^{\frac{p-q}{q}}. \label{eq:39}
\end{align}
Using  doubling condition \eqref{eq: 3.3a} for $\TR$ in the second term $F_{2,\TR}$, we get
\begin{align}
     F_{2,\TR}(x) & = \int_{B^c(0, R)} \tau_{-x}^\alpha |f|(y) \frac{\TR(|y|)}{(1 + |y|)^\gamma}  d\mu_\alpha(y) \nonumber
     \\
     & \leq \int_{B^c(0, R)} \tau_{-x}^\alpha|f|(y)  \frac{\TR(|y|)}{(|y|)^\gamma}  d\mu_\alpha(y) \nonumber
    \\
    & = \sum_{k = 0}^{\infty} \int_{2^k R \leq |y| < 2^{k+1}R} \frac{\TR(|y|)}{(|y|)^\gamma} \tau_{-x}^\alpha |f|(y)  d\mu_\alpha(y) \nonumber
   \\
     & \leq C\sum_{k= 0}^{\infty}   \frac{\TR(2^k R)}{( 2^k R)^\gamma}\int_{2^k R \leq |y| < 2^{k+1}R} \tau_{-x}^\alpha |f|(y)  d\mu_\alpha(y). \nonumber
    \end{align}
Now from \eqref{eq: 4.5b}, we can write
\begin{align}
  F_{2,\TR}(x)  & \leq C ||f||_{L^{p, \phi}(\mathbb{R}, d\mu_\alpha)} \sum_{k= 0}^{\infty} \frac{\TR(2^k R)}{( 2^k R)^\gamma} \phi(2^{k} R) (2^k R)^{d_\alpha} \nonumber
     \\
     & \leq C ||f||_{L^{p, \phi}(\mathbb{R}, d\mu_\alpha)} \sum_{k= 0}^{\infty} \frac{\TR(2^{k+1} R) \phi(2^{k+1} R)}{(2^k R)^{\gamma - d_\alpha }} \nonumber
     \\
     &\leq C ||f||_{L^{p, \phi}(\R, d\mu_\alpha)} \sum_{k=0}^{\infty} \int_{2^k R}^{2^{k+1} R} \frac{\TR(t) \phi(t)}{t^{\gamma - d_\alpha + 1}} dt \nonumber
      \\
     &\leq C ||f||_{L^{p, \phi}(\R, d\mu_\alpha)} \int_{R}^{\infty} \frac{\TR(t) \phi(t)}{t^{\gamma - d_\alpha + 1}} dt, \nonumber
\end{align}
since
\begin{align}
    \int_{2^k R}^{2^{k+1}R} \frac{\TR(t)\phi(t)}{t^{\gamma - d_\alpha + 1}} dt &\geq C\frac{\TR(2^{k+1}R)\phi(2^{k+1}R)}{(2^{k+1}R)^{\gamma - d_\alpha +1}} 2^k R \geq C\frac{\TR(2^{k+1}R)\phi(2^{k+1} R)}{(2^k R)^{\gamma - d_\alpha }} . \nonumber
\end{align}
Now using \eqref{eq:37}, we have
\begin{align}
   F_{2,\TR}(x) \leq C ||f||_{L^{p, \phi}(\R, d\mu_\alpha)} (\phi(R))^{\frac{p}{q}}. \label{eq: 4.7}
\end{align}
Adding equation \eqref{eq:39} and \eqref{eq: 4.7}, we get
\begin{align}
    |I_{\TR, \gamma}^\alpha f(x)| \leq C (M^\alpha f(-x) \phi(R)^{\frac{p-q}{q}} + ||f||_{L^{p, \phi}(\R, d\mu_\alpha)} (\phi(R))^{\frac{p}{q}}). \label{eq: 4.7a}
\end{align}
Note that the above inequality \eqref{eq: 4.7a} holds $\forall~~ R > 0$. Suppose that $f \neq 0$ and $M^\alpha f$ is finite everywhere. Then since $\phi$ is surjective, for each $x \in \R$, we can choose $R_x > 0$ such that $\phi(R_x) = M^\alpha f(-x) ||f||_{L^{p, \phi}(\R, d\mu_\alpha)}^{-1}$. The above inequality \eqref{eq: 4.7a} holds $\forall~~ R >0$ by taking $R = R_x$
\begin{align}
    |I_{\TR, \gamma}^\alpha f(x)| \leq C ||f||_{L^{p, \phi}(\R, d\mu_\alpha)}^{\frac{q-p}{q}}(M^\alpha f(-x))^{\frac{p}{q}} .\nonumber
\end{align}
Now for every $x_0 \in \mathbb{R}$ and for every $r > 0$, consider 
\begin{align}
    &\left(\int_{B(0,r)}\tau_{x_o}^\alpha |I_{\TR, \gamma}^\alpha f|^q(x) d\mu_\alpha(x)\right)^\frac{1}{q} \nonumber\\
    &= \left (\int_{\mathbb{R}} |I_{\TR, \gamma}^\alpha f|^q(x) \tau_{-x_o}^\alpha \chi_{B(0,r)}(x) d\mu_\alpha(x) \right)^\frac{1}{q} \nonumber
    \\
    &\leq C ||f||_{L^{p, \phi}(\R, d\mu_\alpha)}^{\frac{q-p}{q}} \left( \int_{\mathbb{R}}  |M^\alpha f|^p(-x) \tau_{-x_o}^\alpha \chi_{B(0,r)}(x) d\mu_\alpha(x) \right)^\frac{1}{q} \nonumber
    \\
     &\leq C ||f||_{L^{p, \phi}(\R, d\mu_\alpha)}^{\frac{q-p}{q}} \left( \int_{B(0,r)} \tau_{x_o}^\alpha |M^\alpha f|^p(-x) d\mu_\alpha(x) \right)^\frac{1}{q}. \nonumber
\end{align}
Dividing both sides by $(\phi(r))^{\frac{p}{q}} r^{\frac{d_\alpha}{q}}$, we obtain
 \begin{align}
&\frac{\left(\int_{B(0,r)}\tau_{x_o}^\alpha |I_{\TR,\gamma}^\alpha f|^q(x)d\mu_\alpha(x)\right)^{\frac{1}{q}}}{\psi(r) r^{\frac{d_\alpha}{q}} } \nonumber
\\
 &\leq C ||f||_{L^{p, \phi}(\R, d\mu_\alpha)}^{\frac{q-p}{q}} \frac{\left(\int_{B(0, r)}\tau_{x_o}^\alpha |M^\alpha f|^p(-x) d\mu_\alpha(x)\right)^{\frac{1}{q}}}{(\phi(r))^{\frac{p}{q}}r^{\frac{d\alpha}{q}}}, \nonumber
\end{align}
where $\psi(r) = (\phi(r))^{\frac{p}{q}}$. Taking the supremum over $r > 0$ and $x_0 \in \mathbb{R}$, we arrive at
\begin{align}
        ||I_{\TR, \gamma }^\alpha f||_{L^{q, \psi}(\R, d\mu_\alpha)} \leq C ||f||_{L^{p, \phi}(\R, d\mu_\alpha)}^{\frac{q-p}{q}} ||M^\alpha f ||_{L^{p, \phi}(\R, d\mu_\alpha)}^{\frac{p}{q}}. \nonumber
\end{align}
Finally, by using the boundedness of Dunkl-type maximal operator $M^\alpha$ on $L^{p, \phi}(\R, d\mu_\alpha)$ from Theorem \ref{tm: 3.1a}, we obtain the desired result
\begin{align}
    ||I_{\TR, \gamma }^\alpha f||_{L^{q, \psi}(\R, d\mu_\alpha)} \leq C ||f||_{L^{p, \phi}(\R, d\mu_\alpha)}, \nonumber
 \end{align}
 completing the proof.\qed


\section{Boundedness of Generalized  Dunkl-Type Fractional Integral Operator on Generalized Dunkl-type Morrey Space}
\label{Sec:6}
In this section, we prove the  generalized Dunkl-type fractional integral operators $T_{\rho}^\alpha$ is bounded on the generalized Dunkl-type  Morrey space \eqref{eq: 3.2a}.
First, we define the generalized Dunkl-type fractional integral operator $T_{\rho}^\alpha$ by 
\begin{align}
    T_{\rho}^\alpha f(x) := \int_\R  f(y) \tau_x^\alpha \frac{\rho(|y|)}{|y|^{d_\alpha}} d\mu_\alpha(y), \label{eq:38}
\end{align}
where $\rho : \R^{+} \to \R^{+}$ satisfies the doubling condition \eqref{eq: 3.3a} with the condition 
\begin{align}
      \int_0^1 \frac{\rho(t)}{t} dt < \infty.  \label{eq: 39}
\end{align}
Here $K_\rho^\alpha(x) = \frac{\rho(|x|)}{|x|^{d_\alpha}}$ is known as the generalized Dunkl-type Riesz kernel. If $\rho(t) = t^{\beta}$, $0< \beta < d_\alpha$, then we have the Dunkl type fractional integral operator $I_{\beta}^\alpha = T_{\rho}^\alpha$ with Dunkl-type Riesz kernel $K_{\rho}^\alpha = K_{\beta}^\alpha$.
\begin{theorem}
     Let $\rho$ and $\phi$ satisfy the doubling condition \eqref{eq: 3.3a}. Let $\phi$ be surjective and for $1 < p < q < \infty$, it satisfies
    \begin{align}
        \phi(r) \int_0^r \frac{\rho(t)}{t} dt + \int_r^\infty \frac{\rho(t) \phi(t)}{t} dt \leq C (\phi(r))^{\frac{p}{q}}, \label{eq:41}
    \end{align}
    for all $r > 0$. Then the generalized Dunkl-type fractional integral operator is bounded from $L^{p, \phi}(\mathbb{R}, d\mu_\alpha)$ to $L^{q, \psi}(\mathbb{R}, d\mu_\alpha)$ i.e. 
    \begin{align}
        ||T^\alpha_{\rho}f ||_{L^{q, \psi}(\R, d\mu_\alpha)} \leq C ||f||_{L^{p, \phi}(\R, d\mu_\alpha)}, \label{eq: 42}
    \end{align}
    where $\psi(r) = (\phi(r))^{\frac{p}{q}}$.
    \end{theorem}
    \textbf{Proof.} Let $R > 0$ be given. Then for $f \in {L}^{p, \phi}(\R, d\mu_\alpha)$, we write
 \begin{align}
     |T_{\rho}^\alpha f (x)| &\leq  \left( \int_{B(0,R)} + \int_{B^c(0,R)} \right) |f|(y) \tau_{x}^\alpha \frac{\rho(|y|)}{|y|^{d_\alpha}} d\mu_\alpha(y) \nonumber
   \\
   & =  \left( \int_{B(0,R)} + \int_{B^c(0,R)} \right) \tau_{-x}^\alpha |f|(y) \frac{\rho(|y|)}{|y|^{d_\alpha}} d\mu_\alpha(y) = T_{1}(x) + T_{2}(x) , \nonumber
 \end{align} 
where \begin{align}
    T_{1}(x) := \int_{B(0, R)}  \tau_{-x}^\alpha |f|(y)  \frac{\rho(|y|)}{|y|^{d_\alpha}} d\mu_\alpha(y) \nonumber
    \end{align}
    and
    \begin{align}
    T_{2}(x) := \int_{B^c(0, R)}  \tau_{-x}^\alpha |f|(y)  \frac{\rho(|y|)}{|y|^{d_\alpha}}  d\mu_\alpha(y). \nonumber
\end{align}
 Then
\begin{align}
    T_{1}(x) & = \int_{B(0, R)}  \tau_{-x}^\alpha |f|(y) \frac{\rho(|y|)}{|y|^{d_\alpha}}d\mu_\alpha(y) \nonumber
    \\
    & = \sum_{k= -\infty}^{-1} \int_{2^k R \leq |y| < 2^{k+1}R}  \tau_{-x}^\alpha |f|(y)  \frac{\rho(|y|)}{|y|^{d_\alpha}} d\mu_\alpha(y) \nonumber
    \\
    & \leq C \sum_{k= -\infty}^{-1} \frac{\rho(2^k R)}{( 2^k R)^{d_\alpha}}\int_{2^k R \leq |y| < 2^{k+1}R}  \tau_{-x}^\alpha |f|(y)  d\mu_\alpha(y) \nonumber
    \\
    &\leq C M^\alpha f(-x) \sum_{k= -\infty}^{-1} \rho(2^k R). \label{eq: 6.6a}
\end{align}
Since 
\begin{align}
    \int_{2^k R}^{2^{k+1}R} \frac{\rho(t)}{t} dt \geq C \rho(2^k R)  \int_{2^k R}^{2^{k+1}R} \frac{1}{t} dt = C \rho(2^k R) \ln 2, \nonumber
\end{align}
using \eqref{eq:41} and the above estimate in \eqref{eq: 6.6a}, we have
\begin{align}
    T_{1}(x) &\leq C M_\alpha f(-x) \sum_{k=-\infty}^{-1} \int_{2^k R}^{2^{k+1}R} \frac{\rho(t)}{t} dt = C M^\alpha f(-x) \int_0^R \frac{\rho(t)}{t} dt \nonumber
    \\ 
    & \leq C M^\alpha f(-x) (\phi(R))^{\frac{p-q}{q}}. \label{eq:43}
\end{align}
Using  doubling condition \eqref{eq: 3.3a} for $\rho$ in the second term $T_{2}$ , we get
\begin{align}
     T_{2}(x) & = \int_{B^c(0, R)}   \tau_{-x}^\alpha|f|(y) \frac{\rho(|y|)}{|y|^{d_\alpha}}  d\mu_\alpha(y) \nonumber
     \\
    & = \sum_{k = 0}^{\infty} \int_{2^k R \leq |y| < 2^{k+1}R} \frac{\rho(|y|)}{|y|^{d_\alpha}}  \tau_{-x}^\alpha |f|(y)  d\mu_\alpha(y) \nonumber
    \\
    & \leq C\sum_{k= 0}^{\infty} \frac{\rho(2^k R)}{( 2^k R)^{d_\alpha}}\int_{2^k R \leq |y| < 2^{k+1}R}  \tau_{-x}^\alpha |f|(y)  d\mu_\alpha(y). \nonumber
    \end{align}
Now from \eqref{eq: 4.5b}, we can write
\begin{align}
  T_{2}(x)  & \leq C ||f||_{L^{p, \phi}(\mathbb{R}, d\mu_\alpha)}\sum_{k= 0}^{\infty} \frac{\rho(2^k R)}{( 2^k R)^{d_\alpha}} \phi(2^{k+1} R) (2^k R)^{d_\alpha} \nonumber
     \\
      & \leq C ||f||_{L^{p, \phi}(\mathbb{R}, d\mu_\alpha)}\sum_{k= 0}^{\infty} \rho(2^{k+1} R) \phi(2^{k+1} R) \nonumber
     \\
     &\leq C ||f||_{L^{p, \phi}(\R, d\mu_\alpha)} \sum_{k=0}^{\infty} \int_{2^k R}^{2^{k+1} R} \frac{\rho(t) \phi(t)}{t} dt \nonumber
     \\
      &= C ||f||_{L^{p, \phi}(\R, d\mu_\alpha)} \int_{R}^{ \infty} \frac{\rho(t) \phi(t)}{t} dt, \nonumber
\end{align}
since
\begin{align}
    \int_{2^k R}^{2^{k+1}R} \frac{\rho(t)\phi(t)}{t} dt &\geq C \rho(2^{k+1}R)\phi(2^{k+1}R) \int_{2^k R}^{2^{k+1}R} \frac{1}{t} dt \nonumber
    \\
    &\geq C \rho(2^{k+1}R)\phi(2^{k+1} R) \ln 2 .\nonumber
\end{align}
Now using \eqref{eq:41} again, we have
\begin{align}
   T_{2}(x) \leq C ||f||_{L^{p, \phi}(\R, d\mu_\alpha)} (\phi(R))^{\frac{p}{q}}. \label{eq:44}
\end{align}
Adding equations \eqref{eq:43} and \eqref{eq:44}, we get
\begin{align}
    |T_{\rho}^\alpha f(x)| \leq C (M^\alpha f(-x) \phi(R)^{\frac{p-q}{q}} + ||f||_{L^{p, \phi}(\R, d\mu_\alpha)} (\phi(R))^{\frac{p}{q}}). \label{eq: 6.8a}
\end{align}
Note that the above inequality \eqref{eq: 6.8a} holds $\forall~ R > 0$. Suppose that $f \neq 0$ and $M^\alpha f$ is finite everywhere. Then since $\phi$ is surjective, for each $x \in \R$ choose $R_x > 0$ such that $\phi(R_x) = M^\alpha f(-x) ||f||_{L^{p, \phi}(\R, d\mu_\alpha)}^{-1}$. Then, by taking $R = R_x$ in \eqref{eq: 6.8a} we arrive at
\begin{align}
    |T_{\rho}^\alpha f(x)| \leq C ||f||_{L^{p, \phi}(\R, d\mu_\alpha)}^{\frac{q-p}{q}}(M^\alpha f(-x))^{\frac{p}{q}} .\nonumber
\end{align}
Now for every $x_0 \in \mathbb{R}$ and for every $r > 0$, consider
\begin{align}
    &\left(\int_{B(0,r)}\tau_{x_o}^\alpha |T_{\rho}^\alpha f|^q(x) d\mu_\alpha(x)\right)^\frac{1}{q} \nonumber\\ 
    &= \left (\int_{\mathbb{R}} |T_{\rho}^\alpha f|^q(x) \tau_{-x_o}^\alpha \chi_{B(0,r)}(x) d\mu_\alpha(x) \right)^\frac{1}{q} \nonumber
    \\
    &\leq C ||f||_{L^{p, \phi}(\R, d\mu_\alpha)}^{\frac{q-p}{q}} \left( \int_{\mathbb{R}}  |M^\alpha f|^p(-x) \tau_{-x_o}^\alpha \chi_{B(0, r)}d\mu_\alpha(x) \right)^\frac{1}{q} \nonumber
    \\
    &\leq C ||f||_{L^{p, \phi}(\R, d\mu_\alpha)}^{\frac{q-p}{q}} \left( \int_{B(0,r)} \tau_{x_o}^\alpha |M^\alpha f|^p(-x) d\mu_\alpha(x) \right)^\frac{1}{q}. \nonumber
\end{align}
Dividing both sides by $(\phi(r))^{\frac{p}{q}} r^{\frac{d_\alpha}{q}}$, we obtain
 \begin{align}
&\frac{\left(\int_{B(0,r)}\tau_{x_o}^\alpha |T_{\rho}^\alpha f|^q(x)d\mu_\alpha(x)\right)^{\frac{1}{q}}}{\psi(r) r^{\frac{d_\alpha}{q}} } \nonumber
\\
 &\leq C ||f||_{L^{p, \phi}(\R, d\mu_\alpha)}^{\frac{q-p}{q}} \frac{\left(\int_{B(0, r)}\tau_{x_o}^\alpha |M^\alpha f|^p(-x) d\mu_\alpha(x)\right)^{\frac{1}{q}}}{(\phi(r))^{\frac{p}{q}}r^{\frac{d_\alpha}{q}}}, \nonumber
\end{align}
where $ \psi(r) = (\phi(r))^\frac{p}{q}$. Now, we take the supremum over $r > 0$ and $x_0 \in \mathbb{R}$ to get
\begin{align}
        ||T_{\rho}^\alpha f||_{L^{q, \psi}(\R, d\mu_\alpha)} \leq C ||f||_{L^{p, \phi}(\R, d\mu_\alpha)}^{\frac{q-p}{q}} ||M^\alpha f ||_{L^{p, \phi}(\R, d\mu_\alpha)}^{\frac{p}{q}}. \nonumber
\end{align}
Finally, by using the boundedness of Dunkl-type maximal operator $M^\alpha$ on $L^{p, \phi}(\mathbb{R}, d\mu_\alpha)$ from Theorem \ref{tm: 3.1a}, we get 
\begin{align}
 ||T_{\rho}^\alpha f||_{L^{q, \psi}(\R, d\mu_\alpha)} \leq C ||f||_{L^{p, \phi}(\R, d\mu_\alpha)}. \nonumber
 \end{align}
 This completes the proof.\qed


\section{Inequalities for the modified version of the generalized  Dunkl-Type fractional Integral Operator in Dunkl-Type $BMO_{\phi}$ space } \label{Sec:7}
In this section, we prove the boundedness of the modified version of the generalized Dunkl-type fractional integral operator in Dunk-type $BMO_\phi$  space $BMO_\phi(\mathbb{R}, d\mu_\alpha)$.
For a function $\phi : (0, +\infty) \to (0, + \infty)$, we define \\$BMO_{\phi}(\mathbb{R}, d\mu_\alpha)$ to be the space of all locally integrable functions $f$ on $\mathbb{R}$ such that 
\begin{align}
    ||f||_{BMO_{\phi}(\mathbb{R}, d\mu_\alpha)} &:= \sup_{\substack{{r>0}\\{x \in \R}}} \frac{1}{\phi(r)} \frac{1}{\mu_{\alpha}(B(0,r))} \int_{B(0,r)} |\tau_x^\alpha f(y) -f_{B(0,r)}(x) | d\mu_\alpha(y) \nonumber
    \\
    & < \infty, \label{eq: 7.1}
\end{align}
where 
\begin{align}
f_{B(0,r)}(x) = \frac{1}{\mu_\alpha(B(0,r))} \int_{B(0,r)} \tau_x^\alpha f(y) d\mu_\alpha(y). \nonumber
\end{align}
If $\phi(r)\equiv 1$, then $BMO_\phi(\mathbb{R}, d\mu_\alpha)= BMO(\mathbb{R}, d\mu_\alpha)$( see \cite{guliyev2009fractional}). 
Next, for a function $\rho : \mathbb{R}^+ \to \mathbb{R}^+$, we define the modified version of the generalized Dunkl-type fractional integral operator $T_\rho^\alpha$ by
\begin{align}
    \tilde{T}_\rho^\alpha f(x) : = \int_{\mathbb{R}} f(y) \left( \tau_x^\alpha \frac{\rho(|y|)}{|y|^{d_\alpha}} - \frac{\rho(|y|)}{|y|^{d_\alpha}} (1- \chi_{B_0}(y)) \right) d\mu_\alpha(y), \label{eq: 7.2}
\end{align}
where $B_0 = B(0,1)$ denotes the ball center at the origin and of radius $1$ and  $\chi_{B_0}$ is the characteristic function of $B_0 = B (0,1)$. In this definition, we assume that $\rho$ satisfies  \eqref{eq: 3.3a}, \eqref{eq: 39} and the following conditions: 
    \begin{align}
         \int_r^\infty \frac{\rho(t)}{t^2} dt \leq C'_1 \frac{\rho(r)}{r} ~~ for ~~ all~~ r > 0, \label{eq: 45}\\
       \frac{1}{2} \leq \frac{r}{s} \leq 2 \implies \left|\frac{\rho(r)}{r^{d_\alpha}} - \frac{\rho(s)}{s^{d_\alpha}}\right| \leq C''_1 |r- s| \frac{\rho(s)}{s^{d_\alpha + 1}}, \label{eq: 46}
    \end{align}
    where $C'_1, C''_2 > 0$ are independent of $r,s > 0$.
    For example, the function $\rho(r) = r^\beta$ satisfies \eqref{eq: 3.3a}, \eqref{eq: 39} and \eqref{eq: 46} for $0 < \beta < d_\alpha$, and also satisfies \eqref{eq: 45} for $ 0 < \beta < 1. $
    
    Now we state in the following the main theorem of this section.
     \begin{theorem}\label{eq: t8}
      Let $\rho$ satisfy \eqref{eq: 3.3a}, \eqref{eq: 39}, \eqref{eq: 45}, \eqref{eq: 46}. Let $\phi$ and $\psi$ be almost increasing, $\phi(r) \sim \phi(2r) $ and $\psi(r) \sim \psi(2r)$. If, for all $r > 0$,
      \begin{align}
          \int_r^\infty \frac{\rho(t) \phi(t)}{t^2} dt \leq A \frac{\rho(r) \phi(r)}{r}, \label{eq: 47}\\
          \int_0^r \frac{\rho(t)}{t} dt ~\phi(r) \leq A' \psi(r), \label{eq: 48}
      \end{align}
     where $A, A' > 0$ are constants, then $\tilde{T}_\rho^\alpha$ is bounded from $BMO_\phi(\mathbb{R}, d\mu_\alpha)$ to $BMO_\psi(\mathbb{R}, d\mu_\alpha)$.
  \end{theorem}
   
    \begin{remark}
     Since the modified version of the generalized fractional integral operator $T_\rho^\alpha $ \eqref{eq: 7.2} can be written as 
    \begin{align}
        \tilde{T}_\rho^\alpha f(x) = f *_\alpha \frac{\rho(|\cdot|)}{|\cdot|^{d_\alpha}}(x) - f *_\alpha \frac{\rho(|\cdot|)}{|\cdot|^{d_\alpha}}(1 - 
  \chi_{B_0}) (0), \nonumber 
    \end{align}
    hence by using Lemma \ref{7.1}, we can write 
    \begin{align}
        \tau_{x_0}^\alpha \tilde{T}_\rho^\alpha f(x) &= \tau_{x_0}^\alpha f *_\alpha \frac{\rho(|\cdot|)}{|\cdot|^{d_\alpha}} (x) - \tau_{x_0}^\alpha f *_\alpha \frac{\rho(|\cdot|)}{|\cdot|^{d_\alpha}}(1 - \chi_{B_0}) (0) \nonumber
        \\
        &= \int_{\mathbb{R}} \tau_{x_0}^\alpha f(y) \left( \tau_x^\alpha \frac{\rho(|y|)}{|y|^{d_\alpha}} - \frac{\rho(|y|)}{|y|^{d_\alpha}}(1 - \chi_{B_0}(y)) \right) d\mu_\alpha(y). \label{eq: 7a}
    \end{align} 
     \end{remark}
      In order to prove the above Theorem \ref{eq: t8}, we first prove the following lemmas:
     \begin{lemma} \label{eq: l7.1}
        For $2|x| \leq |y|$, the following inequalities are valid:
        \begin{itemize}
         \item [(i)] $
            \left|\tau_x^\alpha \frac{\rho(|y|)}{|y|^{d_\alpha}} - \frac{\rho(|y|)}{|y|^{d_\alpha}} \right| \leq C |x| \frac{\rho(|y|)}{|y|^{d_\alpha + 1}}, 
            $\\
            \item [(ii)] $
            \left|\tau_x^\alpha \frac{\rho(|y|)}{|y|^{d_\alpha}} \right| \leq C \frac{\rho(|y|)}{|y|^{d_\alpha}}, 
             $
         \end{itemize}
        where $C > 0$ is constant and $\rho$ satisfies the doubling condition \eqref{eq: 3.3a}.
    \end{lemma}
  \textbf{Proof.}  Since
  \begin{align}
  \int_0^\pi (1- cos \theta) sin^{2\alpha} \theta d\theta & = 2^{2 \alpha + 1} \int_0^\pi sin^{2 \alpha + 2} \frac{\theta}{2} cos^{2 \alpha} \frac{\theta}{2} d\theta \nonumber
  \\
  & = 2^{2 \alpha +1} \frac{\Gamma(\alpha + \frac{1}{2}) \Gamma( \alpha + \frac{3}{2} )}{\Gamma(2 \alpha + 2)} \nonumber
  \\
  & = \frac{\sqrt{\pi} \Gamma(\alpha + \frac{1}{2})}{\Gamma(\alpha + 1)} = \frac{1}{c_\alpha}, \nonumber
  \end{align}
we get from \eqref{2.1b}
    \begin{align}
        \left| \tau_x^\alpha \frac{\rho(|y|)}{|y|^{d_\alpha}} - \frac{\rho(|y|)}{|y|^{d_\alpha}} \right|
        & =\left| c_\alpha \int_0^\pi \frac{\rho(|(x,y)_\theta|)}{|(x,y)_\theta|^{d_\alpha}} (1 - cos \theta) sin^{2\alpha} \theta d\theta - \frac{\rho(|y|)}{|y|^{d_\alpha}}\right| \nonumber
        \\
        & = c_\alpha \left| \int_0^\pi \frac{\rho(|(x,y)_\theta|)}{|(x,y)_\theta|^{d_\alpha}} (1 - cos \theta) sin^{2\alpha} \theta d\theta - \frac{1}{c_\alpha}\frac{\rho(|y|)}{|y|^{d_\alpha}}\right| \nonumber\\
        & \leq 2 c_\alpha \int_0^\pi \left| \frac{\rho(|(x,y)_\theta|)}{|(x,y)_\theta|^{d_\alpha}} - \frac{\rho(|y|)}{|y|^{d_\alpha}} \right | sin^{2\alpha}\theta d\theta. \label{eq: 2c}
    \end{align}
   Now using the assumption $2 |x| \leq |y|$, it is easy to see 
    \begin{align}
        |(x,y)_\theta| \leq |x| + |y| \leq \frac{3}{2} |y| , \hspace{1.5cm} |(x,y)_\theta| \geq |y| - |x| \geq \frac{1}{2} |y|. \nonumber
    \end{align}
    This implies
    \begin{align}
        |(x,y)_\theta| - |y| \leq |x|, \hspace{1.5cm} |y| - |(x,y)_\theta| \leq |x|. \nonumber
    \end{align}
    Hence, by combining these inequalities, we obtain 
    \begin{align}\frac{1}{2}|y| \leq |(x,y)_\theta| \leq \frac{3}{2} |y| \hspace{1cm} and \hspace{1cm} ||(x,y)_\theta| - |y|| \leq |x|. \nonumber
    \end{align}
    \\
  Now taking $ r=|(x,y)_\theta| $ and $ s= |y|$ in \eqref{eq: 46}, we get from \eqref{eq: 2c}
  \begin{align}
       \left| \tau_x^\alpha \frac{\rho(|y|)}{|y|^{d_\alpha}} - \frac{\rho(|y|)}{|y|^{d_\alpha}} \right|  &\leq 2 c_\alpha \int_0^\pi \left| |(x,y)_\theta| - |y|  \right | \frac{\rho(|y|)}{|y|^{d_\alpha + 1}} sin^{2\alpha}\theta d\theta \nonumber
       \\
       & \leq  2 C''_1 c_\alpha |x| \frac{\rho(|y|)}{|y|^{d_\alpha +1}} \int_0^\pi sin^{2 \alpha} \theta d\theta = C |x| \frac{\rho(|y|)}{|y|^{d_\alpha + 1}}. \nonumber
  \end{align}
  Hence $(i)$ is proved. The inequality $(ii)$ follows similarly.
  \qed   

  \begin{lemma} \label{eq: l7.2}
      Under the assumptions in Theorem \ref{eq: t8}, there exists a constant $C > 0$ such that for all $x \in \mathbb{R}$ and for all $r>0$,
      \begin{align}
          \int_{B^c(0,r)} \frac{\rho(|y|)}{|y|^{d_\alpha + 1}} |\tau_x^\alpha f(y) - f_{B(0,r)}(x)| d\mu_\alpha(y) \leq C \frac{\rho(r) \phi(r)}{r} ||f||_{BMO_\phi(\mathbb{R}, d\mu_\alpha)}. \nonumber
      \end{align}
  \end{lemma}
  \textbf{Proof.} Let $k \in \mathbb{N}$. Then using the assumption on $\phi$, we get
  \begin{align}
      &|f_{B(0, 2^k r)}(x) - f_{B(0, 2^{k+1} r)}(x)| \nonumber
      \\
      &\leq \frac{1}{\mu_\alpha(B(0, 2^k r))} \int_{B(0, 2^k r)} |\tau_x^\alpha f(y) - f_{B(0, 2^{k+1}r)}(x)| d\mu_\alpha(y) \nonumber\\
      &\leq \frac{1}{\mu_\alpha(B(0, 2^k r))} \int_{B(0, 2^{k+1} r)} |\tau_x^\alpha f(y) - f_{B(0, 2^{k+1}r)}(x)| d\mu_\alpha(y) \nonumber\\
      & \leq 2^{d_\alpha} \phi(2^{k+1}r) ||f||_{BMO_\phi(\mathbb{R}, d\mu_\alpha)} \nonumber \\
      & \leq C \int_{2^k r}^{2^{k+1}r} \frac{\phi(s)}{s} ds ||f||_{BMO_\phi(\mathbb{R}, d\mu_\alpha)}. \nonumber 
  \end{align}
  Hence for all $j \geq 1$, by taking $k =0, 1, 2, ..., j-1$ in the above estimate, we get
  \begin{align}
      &\frac{1}{\mu_\alpha(B(0, 2^j r))} \int_{B(0, 2^j r)} |\tau_x^\alpha f(y) - f_{B(0, r)}(x)| d\mu_\alpha(y) \nonumber \\
      &\leq \frac{1}{\mu_\alpha(B(0, 2^j r))} \int_{B(0, 2^j r)} |\tau_x^\alpha f(y) - f_{B(0, 2^jr)}(x)| d\mu_\alpha(y) \nonumber 
      \\
      &\hspace{4.5 cm}+ |f_{B(0,r)}(x) - f_{B(0, 2^j r)}(x)|\nonumber 
      \\
      &\leq \phi(2^j r) ||f||_{BMO_{\phi}(\mathbb{R}, d\mu_\alpha)} + C \int_{r}^{2^j r} \frac{\phi(s)}{s} ds ||f||_{BMO_\phi(\mathbb{R}, d\mu_\alpha)} \nonumber
      \\
      &\leq C \int_r^{2^j r} \frac{\phi(s)}{s} ds ||f||_{BMO_\phi(\mathbb{R}, d\mu_\alpha)}, \label{eq: 7.8c}
  \end{align}
  since
  \begin{align*}
      C \int_r^{2^j r} \frac{\phi(s)}{s} ds \geq C \int_{2^{j-1} r}^{2^j r} \frac{\phi(s)}{s} ds \geq \phi(2^j r).
  \end{align*}
 Now using \eqref{eq: 45} and \eqref{eq: 47}, we have 
  \begin{align}
    & \int_{B^c(0, r)} \frac{\rho(|y|)}{|y|^{d_\alpha + 1}} |\tau_x^\alpha f(y) - f_{B(0,r)}(x)| d\mu_\alpha(y) \nonumber\\
    & = \sum_{j=1}^\infty \int_{2^{j-1}r \leq |y| < 2^j r} \frac{\rho(|y|)}{|y|^{d_\alpha +1}} |\tau_x^\alpha f(y) - f_{B(0,r)}(x)| d\mu_\alpha(y) \nonumber \\
    & \leq C \sum_{j=1}^\infty \frac{\rho(2^j r)}{(2^j r)^{d_\alpha + 1}} \int_{B(0, 2^j r)} |\tau_x^\alpha f(y) - f_{B(0,r)}(x)| d\mu_\alpha(y). \nonumber
    \end{align}
    Since 
     \begin{align*}
      \int_{2^{j-1}r}^{2^j r} \frac{\rho(t)}{t^2} \left( \int_r^{2t} \frac{\phi(s)}{s} ds\right) dt \geq C\frac{\rho(2^j r)}{(2^j r)^2} \left( \int_r^{2^{j}r} \frac{\phi(s)}{s} ds \right) 2^{j-1} r
  \end{align*} 
  and by \eqref{eq: 7.8c}, we get
    \begin{align}
      &\int_{B^c(0, r)} \frac{\rho(|y|)}{|y|^{d_\alpha + 1}} |\tau_x^\alpha f(y) - f_{B(0,r)}(x)| d\mu_\alpha(y) \nonumber
      \\
      & \leq C \sum_{j=1}^\infty \frac{\rho(2^j r)}{2^j r} \int_r^{2^j r} \frac{\phi(s)}{s} ds ||f||_{BMO_\phi(\mathbb{R}, d\mu_\alpha)} \nonumber
    \\
    &\leq C \int_r^\infty \frac{\rho(t)}{t^2} \left( \int_r^{2t} \frac{\phi(s)}{s} ds \right) dt ||f||_{BMO_\phi(\mathbb{R}, d\mu_\alpha)} \nonumber\\
    & = C\int_r^\infty \left( \int_{s/2}^{\infty} \frac{\rho(t)}{t^2} dt \right) \frac{\phi(s)}{s} ds ||f||_{BMO_\phi(\mathbb{R}, d\mu_\alpha)} \nonumber\\ 
    & \leq C \int_r^\infty \frac{\rho(s)}{s} \frac{\phi(s)}{s} ds ||f||_{BMO_\phi(\mathbb{R}, d\mu_\alpha)} \nonumber
    \\
    &\leq C \frac{\rho(r) \phi(r)}{r} ||f||_{BMO_\phi(\mathbb{R}, d\mu_\alpha)},\nonumber 
  \end{align} 
 which completes the proof of Lemma \ref{eq: l7.2}.
  \qed ~\\
  \textbf{Proof of Theorem \ref{eq: t8}.} Let $f \in BMO_\phi(\mathbb{R}, d\mu_\alpha)$. For given $r >0$, let $\tilde{B}=B(0,2r)$ and suppose $x \in B(0,r)$.
  Consider
  \begin{align}
     & E_{B(0,r)}(x) = \int_\mathbb{R} ( \tau_{x_0}^\alpha f(y) - f_{\tilde{B}}(x_0) ) \left( \tau_x^\alpha \frac{\rho(|y|)}{|y|^{d_\alpha}} - \frac{\rho(|y|)}{|y|^{d_\alpha}} (1 - \chi_{\tilde{B}}(y))\right) d\mu_\alpha(y),  \nonumber
     \\
      & C_{B^1(0,r)}(x_0) = \int_\mathbb{R} ( \tau_{x_0}^\alpha f(y) - f_{\tilde{B}}(x_0) ) \nonumber
      \\
      &\hspace{3.5cm} \times\left( \frac{\rho(|y|)}{|y|^{d_\alpha}}(1 - \chi_{\tilde{B}}(y)) - \frac{\rho(|y|)}{|y|^{d_\alpha}} (1 - \chi_{B_0}(y))\right) d\mu_\alpha(y), \nonumber
      \\
       & C_{B^2(0,r)}(x_0) = \int_\mathbb{R} f_{\tilde{B}}(x_0)  \left( \tau_{x}^\alpha \frac{\rho(|y|)}{|y|^{d_\alpha}}  - \frac{\rho(|y|)}{|y|^{d_\alpha}} (1 - \chi_{B_0}(y))\right) d\mu_\alpha(y), \nonumber
       \\
         & E_{B^1(0,r)}(x) = \int_{\tilde{B}} ( \tau_{x_0}^\alpha f(y) - f_{\tilde{B}}(x_0) ) \tau_x^\alpha \frac{\rho(|y|)}{|y|^{d_\alpha}} d\mu_\alpha(y), \hspace{6.25 cm} \nonumber
         \\
          & E_{B^2(0,r)}(x) = \int_{\tilde{B}^c} ( \tau_{x_0}^\alpha f(y) - f_{\tilde{B}}(x_0) ) \left( \tau_x^\alpha \frac{\rho(|y|)}{|y|^{d_\alpha}} - \frac{\rho(|y|)}{|y|^{d_\alpha}}\right) d\mu_\alpha(y),  \nonumber 
  \end{align}
  where $x_0 \in \mathbb{R}$. Then, from \eqref{eq: 7a}
  \begin{align}
      \tau_{x_0}^\alpha \tilde{T}_\rho^\alpha f(x) - (C_{B^1(0,r)}(x_0) + C_{B^2(0,r)}(x_0)) &= E_{B(0,r)}(x) \nonumber
      \\
      &= E_{B^1(0,r)}(x) + E_{B^2(0,r)}(x). \nonumber
  \end{align}
  Since 
  \begin{align}
     &\left|\frac{\rho(|y|) (1 - \chi_{\tilde{B}}(y))}{|y|^{d_\alpha}} - \frac{\rho(|y|)(1 - \chi_{B_0}(y))}{|y|^{d_\alpha}}\right| \nonumber\\
     & \hspace{3cm}\leq  \begin{cases}
            0, &  |y| < min (1,2r) ~or ~|y| \geq max(1, 2r);\\
            \frac{\rho(|y|)}{|y|^{d_\alpha}}
            , & \text{otherwise},
        \end{cases}\nonumber
          \end{align}
    this implies $C_{B^1(0,r)}(x_0)$ is finite.

Now let us show that $C_{B^2(0,r)}(x_0)$ is finite. For this first we shall show that the following integral is finite:
\begin{align}
    &\int_{\mathbb{R}} \left(\tau_x^\alpha \frac{\rho(|y|)}{|y|^{d_\alpha}} - \frac{\rho(|y|)(1- \chi_{B_0}(y))}{|y|^{d_\alpha}} \right) d\mu_\alpha(y) \nonumber\\
    &= \int_{\mathbb{R}} \left(\tau_x^\alpha \frac{\rho(|y|)}{|y|^{d_\alpha}}- \frac{\rho(|y|)}{|y|^{d_\alpha}}\right) d\mu_\alpha(y) + \int_{B_0} \frac{\rho(|y|)}{|y|^{d_\alpha}} d\mu_\alpha(y). \label{7.10a}
\end{align}
Let us denote $A := \int_{\mathbb{R}} \left(\tau_x^\alpha \frac{\rho(|y|)}{|y|^{d_\alpha}}- \frac{\rho(|y|)}{|y|^{d_\alpha}}\right) d\mu_\alpha(y)$. For sufficiently large $R > 0$, we write $A$ in the form $A = A_1 + A_2 + A_3$, where 
\begin{align}
    &A_1 := \int_{B(0,R)} \tau_x^\alpha \frac{\rho(|y|)}{|y|^{d_\alpha}} d\mu_\alpha(y) -  \int_{B(-x,R)} \frac{\rho(|y|)}{|y|^{d_\alpha}} d\mu_\alpha(y), \nonumber\\
    &A_2 :=  \int_{B(0,R+r)\backslash B(0, R)} \tau_x^\alpha \frac{\rho(|y|)}{|y|^{d_\alpha}} d\mu_\alpha(y) -  \int_{B(0,R+r)\backslash B(-x, R)} \frac{\rho(|y|)}{|y|^{d_\alpha}} d\mu_\alpha(y), \nonumber\\
    &A_3 :=\int_{B^c(0,R+r)} \left( \tau_x^\alpha \frac{\rho(|y|)}{|y|^{d_\alpha}} - \frac{\rho(|y|)}{|y|^{d_\alpha}}\right) d\mu_\alpha(y). \hspace{1.25cm}\nonumber
\end{align}
First we consider the integral $A_3$. In this case $|y| \geq R +r$. Since $R$ is sufficiently large, we can choose $R > r$. Now since $x \in B(0,r)$, $|y| \geq 2|x|$. Hence by Lemma \ref{eq: l7.1}, we have
\begin{align}
    |A_3| &\leq \int_{B^c(0, R+r)} \left| \tau_x^\alpha \frac{\rho(|y|)}{|y|^{d_\alpha}} - \frac{\rho(|y|)}{|y|^{d_\alpha}}\right| d\mu_\alpha(y) \nonumber \\
          & \leq C r \int_{B^c(0, R+r)} \frac{\rho(|y|)}{|y|^{d_\alpha + 1 }} d\mu_\alpha(y) \nonumber\\
          & = C r \int_{R+r}^\infty \frac{\rho(t)}{t^2} dt. \nonumber
\end{align}
The inequality \eqref{eq: 45} implies that the last integral is finite and $|A_3| \to 0$ as $R \to +\infty$. Now we consider the integral $A_2$.  Since $B(0, R+r) \backslash B(0, R) \subseteq B(0, R+ r) \backslash B(0, R -r)$ and $B(0, R +r) \backslash B(-x, R) \subseteq B(0, R +r) \backslash B(0, R -r)$, 
\begin{align}
   |A_2| \leq \int_{B(0, R +r)\backslash B(0,R - r)} \left(\tau_x \frac{\rho(|y|)}{|y|^{d_\alpha}} + \frac{\rho(|y|)}{|y|^{d_\alpha}} \right) d\mu_\alpha(y). \nonumber 
\end{align}
Further since $2|x| \leq 2r \leq R$ for large $R$,  using Lemma \ref{eq: l7.1}, we get
\begin{align}
    |A_2| &\leq C \int_{B(0, R+r)\backslash B(0, R - r)} \frac{\rho(|y|)}{|y|^{d_\alpha}} d\mu_\alpha(y) \nonumber\\
          & \sim ((R+r)^{d_\alpha} - (R - r)^{d_\alpha}) \frac{\rho(R)}{R^{d_\alpha}} \leq C r \frac{\rho(R)}{R}.\nonumber
\end{align}
The last term goes to zero as $R \to +\infty$ i.e. $|A_2| \to 0$ as $R \to +\infty$. 
Now for $A_1$, we have
\begin{align}
    |A_1| &\leq \int_{B(-x, R)} \frac{\rho(|y|)}{|y|^{d_\alpha}} \left| \tau_{-x}^\alpha \chi_{B(0,R)}(y) -1 \right| d\mu_\alpha(y) \nonumber\\
           &\leq 2 \int_{B(-x, R)} \frac{\rho(|y|)}{|y|^{d_\alpha}} d\mu_\alpha(y) \leq \int_{B(0, r+ R)} \frac{\rho(|y|)}{|y|^{d_\alpha}} d\mu_\alpha(y). \nonumber
\end{align}
Thus using  \eqref{eq: 39}, the right-hand side of the above integral is finite, and hence $A_1$ is finite. Hence using again \eqref{eq: 39}, it follows from \eqref{7.10a} that $C_{B^2(0,r)}(x_0)$ is finite.
Now consider
\begin{align}
    &\int_{B(0, 2r)}\left( \int_{B(0,r)} |\tau_{x_0}^\alpha f(y) - f_{B(0,2r)}(x_0)| \tau_x^\alpha \frac{\rho(|y|)}{|y|^{d_\alpha}} d\mu_\alpha(x) \right ) d\mu_\alpha(y) \nonumber\\
    & \leq \int_{B(0, 2r)} |\tau_{x_0}^\alpha f(y) - f_{B(0,2r)}(x_0)| \left( \int_{B(0, r)} \tau_y^\alpha \frac{\rho(|x|)}{|x|^{d_\alpha}} d\mu_\alpha(x) \right) d\mu_\alpha(y) \nonumber \\
     & \leq \int_{B(0, 2r)} |\tau_{x_0}^\alpha f(y) - f_{B(0,2r)}(x_0)| \nonumber
     \\
     & \hspace{3cm} \times \left( \int_{B(-y, r)} \frac{\rho(|x|)}{|x|^{d_\alpha}} \tau_{-y}^\alpha \chi_{B(0,r)}(x) d\mu_\alpha(x) \right) d\mu_\alpha(y) \nonumber \\
     & \leq \int_{B(0, 2r)} |\tau_{x_0}^\alpha f(y) - f_{B(0,2r)}(x_0)| \left( \int_{B(0, 3r)} \frac{\rho(|x|)}{|x|^{d_\alpha}} d\mu_\alpha(x) \right) d\mu_\alpha(y) \nonumber \\
    & = \int_{B(0, 2r)} |\tau_{x_0}^\alpha f(y) - f_{B(0, 2r)}(x_0) | \left( \int_0^{3r} \frac{\rho(t)}{t} dt \right) d\mu_\alpha(y) \nonumber \\
    & \leq C ||f||_{BMO_{\phi}(\mathbb{R}, d\mu_\alpha)} r^{d_\alpha} \phi(r) \int_0^r \frac{\rho(t)}{t} dt .\nonumber
\end{align}
By using \eqref{eq: 48} and Fubini's theorem, it follows that $E_{B^1(0,r)}$ is finite and 
\begin{align}
    \int_{B(0, r)} |E_{B^1(0,r)}(x)| d\mu_\alpha(x) \leq C \psi(r) r^{d_\alpha} ||f||_{BMO_\phi(\mathbb{R}, d\mu_\alpha)}. \label{eq: 49}
\end{align}
For $E_{B^2(0,r)}$, we have
\begin{align}
    |E_{B^2(0,r)}(x)| \leq \int_{B^c(0,2r)} |\tau_{x_0}^\alpha f(y) - f_{B(0,2r)}(x_0)| \left| \tau_x^\alpha \frac{\rho(|y|)}{|y|^{d_\alpha}} - \frac{\rho(|y|)}{|y|^{d_\alpha}} \right| d\mu_\alpha(y). \nonumber
\end{align}
Then Lemma \ref{eq: l7.1} and Lemma \ref{eq: l7.2} imply that $E_{B^2(0,r)}$ is finite and 
\begin{align}
  |E_{B^2(0,r)}(x)| &\leq C r \int_{B^c(0,2r)} |\tau_{x_0}^\alpha f(y) - f_{B(0, 2r)}(x_0)| \frac{\rho(|y|)}{|y|^{d_\alpha +1}} d\mu_\alpha(y)  \nonumber\\
  & \leq C \rho(r) \phi(r) ||f||_{BMO_\phi(\mathbb{R}, d\mu_\alpha)} \nonumber\\
  & \leq C \int_0^r \frac{\rho(t)}{t} dt \phi(r) ||f||_{BMO_\phi(\mathbb{R}, d\mu_\alpha)} \nonumber\\
  & \leq C \psi(r) ||f||_{BMO_\phi(\mathbb{R} , d\mu_\alpha)}. \nonumber
\end{align}
This implies
\begin{align}
    \int_{B(0,r)} |E_{B^2(0,r)}(x)| d\mu_\alpha(x) \leq C \psi(r) r^{d_\alpha}||f||_{BMO_\phi(\mathbb{R} , d\mu_\alpha)}.  \label{eq: 50}
\end{align}
Finally, by \eqref{eq: 49} and \eqref{eq: 50}, we have 
\begin{align}
    \frac{1}{\mu_\alpha(B(0,r))} \int_{B(0,r)} |\tau_{x_0}^\alpha \tilde{T}_\rho^\alpha f(x) - (C_{B^1(0,r)}(x_0) &+ C_{B^2(0,r)}(x_0))| d\mu_\alpha(x)  \nonumber
    \\
    &\leq C ||f||_{BMO_\phi(\mathbb{R}, d\mu_\alpha)} \psi(r), \nonumber
\end{align} 
which implies
\begin{align}
    ||\tilde{T}_\rho^\alpha f ||_{BMO_\psi(\mathbb{R}, d\mu_\alpha)}   \leq C ||f||_{BMO_\phi(\mathbb{R}, d\mu_\alpha)}, \nonumber
\end{align}
completing the proof.\qed

  \bibliographystyle{amsplain}
  \bibliographystyle{abbrv}

\begin{thebibliography}{99}
\bibitem{abdelkefi2007characterization} {Abdelkefi, C.,  Sifi, M.}: Characterization of Besov spaces for the Dunkl operator on the real line. JIPAM. J. Inequal. Pure Appl. Math. 8(3)(2007) 1-21, Article 73.
\bibitem {abdelkefi2007dunkl} {Abdelkefi, C.,  Sifi, M.}:   Dunkl translation and uncentered maximal operator on the real line. Int. J. Math. Math. Sci. vol. 2007, Article ID 87808, 9 pages, 2007.
\bibitem{adams1975note} {Adams, D.R.}: A note on Riesz potentials. Duke Math. J. 42 (1975), pp. 765-778.
\bibitem{sps}{Adhikari, S., Parui, S.}: The boundedness of Dunkl Bessel Riesz operators on Dunkl-type Morrey spaces. Indian Journal of Pure and Applied Mathematics (2025): 1-14.
\bibitem {amri2012riesz} {Amri, B., Sifi, M.}: Riesz transforms for Dunkl transform. Ann. Math. Blaise Pascal. 19 (2012), pp. 247-262. 
 \bibitem {chiarenza1987morrey} {Chiarenza, F., Frasca, M.}: Morrey spaces and Hardy-Littlewood maximal function. Rend. Mat. Appl., VII. Ser. 7 (1987), pp. 679–693.
 \bibitem {de1993dunkl} {De Jeu, M.F.E.}: The Dunkl transform. Invent. Math. 113 ( 1993), pp. 147–162.
\bibitem {dunkl1989differential} {Dunkl, C.F.}: Differential-difference operators associated to reflection groups. Trans. Amer. Math. Soc. 311 (1989), pp. 167-183.
\bibitem {dunkl1992hankel} {Dunkl, C.F.}: Hankel Transforms Associated to Finite Reflection Groups. Hypergeometric Functions on Domains of
Positivity, Jack polynomials, and applications (Tampa, Fla, 1991), of Contemporary Mathematics, Vol. 138, American
Mathematical Society, Providence, RI, USA, pp. 123–138.
\bibitem{dunkl1991integral} {Dunkl, C.F.}: Integral kernels with reflection group invariance. Can. J. Math. 1991 Dec;43(6):1213-27.
\bibitem{eridani2002boundedness} {Eridani}: On the boundedness of a generalized fractional integral on generalized Morrey spaces. Tamkang J. Math. 33, no. 4 (2002): 335-340.
\bibitem{eridani2004generalized} {Eridani, Gunawan, H., Nakai, E.}: On generalized fractional integral operators. Sci. Math. Jpn. 60, no. 3 (2004): 539-550.
\bibitem {dg} {Gorbachev, D.V., Ivanov, V.I., Tikhonov, S. Yu.}: Riesz Potential and maximal function for Dunkl transform. Potential Anal. (2020), https://doi.org/10.1007/s11118-020-09867-z.
\bibitem {guliyev2010necessary}  {Guliyev, E., Eroglu, A., Mammadov, Y.Y.}: Necessary and Sufficient conditions for the boundedness of Dunkl-type Fractional maximal Operator in the Dunkl-Type Morrey spaces. Abstr. Appl. Anal. Volume 2010 (2010), Art. ID 976493, 10 pp.; doi:10.1155/2010/976493.
\bibitem{Gul}{Guliyev, V.S.,}: Integral operators on function spaces on the homogeneous groups and in domain $\mathbb{G}$. Doctor's Degree Dissertation, Moscow, Mat. Inst. Steklov (1994), 1-329
\bibitem {guliyev2010p} {Guliyev, V.S., Mammadov, Y.Y.}: $(L_p, L_q)$ boundedness of the fractional maximal operator
associated with the Dunkl operator on the real line”, Integral Transforms Spec. Funct. (2010) pp. 1–11.
\bibitem{guliyev2009fractional} {Guliyev, V.S., Mammadov, Y.Y.}: On fractional maximal function and fractional integrals associated with the Dunkl operator on the real line. J. Math. Anal. Appl. (2009) 353(1), 449-459.
\bibitem {hardy1928some} {Hardy, G.H., Littlewood, J.E.}: Some properties of fractional integrals I. Math. Z. 27 (1928), pp. 565–606. 
\bibitem{hardy1932some}{Hardy, G.H., Littlewood, J.E.}: Some properties of fractional integrals II. Math. Z. 34 (1932), pp. 403–439.
\bibitem{2017boundedness} {Idris, M., Gunawan, H.}: The boundedness of generalized Bessel-Riesz operators on generalized Morrey spaces. In Journal of Physics: Conference Series $2017 Oct 1 (Vol. 893, No. 1, p. 012014)$. IOP Publishing.
\bibitem{idris2016boundedness} {Idris, M., Gunawan, H., Eridani}: The boundedness of Bessel-Riesz operators on generalized Morrey spaces. Aust. J. Math. Anal. Appl. 13, No 1(2016), 1-10
\bibitem{lk} {Kamoun, L.}: Besov-type spaces for the Dunkl operator on the real line. J. Comput. Appl. Math. 199, no. 1 (2007): 56-67.
\bibitem{yy} {Mammadov, Y.Y.,  Hasanli, S.A.:} On the boundedness of commutators Dunkl-type maximal operator in the Dunkl-type Morrey spaces. Caspian
J. of Appl. Mathematics, Ecology and Economics, 6, No 1 (2018), 37-52.
\bibitem{mourou2001transmutation} {Mourou, M.A.}: Transmutation operators associated with a Dunkl type differential-difference operator on the real line and certain of their applications. Integral Transform Spec. Funct. 2001 Aug 1;12(1):77-88.
\bibitem{Nagacy}{Nagacy, P., Feuto, J.}: Maximal operator in Dunkl-Fofana spaces. Advances in Pure and Applied Mathematics.(2021). 12. 30-52. 10.21494/ISTE.OP.2021.0647. 
\bibitem{nakai2001generalized} {Nakai, E.}: On generalized fractional integrals. Taiwan. J. Math. 2001 Sep;5(3):587-602.
\bibitem {peetre1969theory} {Peetre, J.}: On the theory of $\mathcal{L}_{p,\lambda}$ spaces. J. Funct. Anal. 4 (1969), pp. 71–87.
\bibitem {rosler1994bessel}{R$\ddot{o}$sler, M.}: Bessel-type signed hypergroups on $\mathbb{R}$. In Probability Measures on Groups and Related Structures, XI (Oberwolfach, 1994), pages 292–304, World Scientific, River edge, NJ, USA, 1995.
\bibitem {rsl2} {R$\ddot{o}$sler, M.}: Dunkl operators: theory and applications", In Orthogonal polynomials and special functions
(Leuven, 2002), volume 1817 of Lecture Notes in Math., pages 93–135. Springer, Berlin, 2003.
\bibitem{Ruz}{Ruzhansky, M., Suragan, D., Yessirkegenov, N.}: Hardy-Littlewood, Bessel-Riesz, and fractional integral operators in anisotropic Morrey and Campanato spaces. Fract. Calc. Appl. Anal. 2018 Jun 26;21(3):577-612. https://doi.org/10.1515/fca-2018-0032.
\bibitem{sifi2002generalized} {Sifi, M., Soltani, F.}: Generalized Fock spaces and Weyl relations for the Dunkl kernel on the real line. J. Math. Anal. Appl. 270, no. 1 (2002): 92-106.
\bibitem {sobolev1938theorem} {Sobolev, S.L.}: On a theorem of functional analysis. Am. Math. Soc., Transl., II. Ser. 34 (1963), 39–68; Translated from Mat. Sb., N. Ser. 4 (1938), 471–497.
\bibitem {soltani2004lp} {Soltani, F.}: $L^p$-Fourier multipliers for the Dunkl operator on the real line. J. Funct. Anal.  209 (1) (2004), pp. 16–35.
\bibitem{solatani2005} {Soltani, F.}: Littlewood–Paley operators associated with the Dunkl operator on $\R$. J. Funct. Anal.  221.1 (2005): 205-225.
\bibitem {sy} {Thangavelu, S., Xu, Y.}: Convolution operator and maximal function for the Dunkl transform. J. Anal. Math. 97 (1)(2005), pp 25-55.
\bibitem {thangavelu2007riesz} {Thangavelu, S., Xu, Y.}: Riesz transform and Riesz potentials for Dunkl transform.  J. Comput. Appl. Math. 199 (2007), pp. 181-195. 
\bibitem{trime2002paley} {Trime'Che, K.}: Paley-Wiener theorems for the Dunkl transform and Dunkl translation operators. Integral Transform Spec. Funct. 2002 Jan 1;13:17-38.
\bibitem{wei} {Wei, M.}: Fractional integral operator and its commutator on generalized Morrey spaces associated with ball Banach function spaces. Fract Calc Appl Anal 26, 2318–2336 (2023). https://doi.org/10.1007/s13540-023-00188-7

\end{thebibliography}
\section*{Acknowledgements}
First author is supported by the University Grants Commission(UGC) with Fellowship No. 211610060698/(CSIRNETJUNE2021). Second author is supported by  seed grant IITM/SG/SWA/94 from Indian Institute of Technology Mandi, India.

\textbf {Data availability statement:} Not applicable

\end{document}